\documentclass{article}

\usepackage[utf8]{inputenc} 
\usepackage[T1]{fontenc}    
\usepackage{hyperref}       
\usepackage{url}            
\usepackage{booktabs}       
\usepackage{amsfonts}       
\usepackage{nicefrac}       
\usepackage{microtype}      
\usepackage{amsmath,amssymb}
\usepackage{mathtools}
\usepackage{comment}
\usepackage{graphicx}
\usepackage{subfigure}
\usepackage{geometry}
\usepackage{bm}

\usepackage{color}

\usepackage[ruled,vlined,linesnumbered]{algorithm2e}
\newtheorem{theorem}{Theorem}

\newtheorem{assumption}{Assumption}
\newtheorem{corollary}{Corollary}
\newtheorem{innercustomthm}{Theorem}
\newenvironment{customthm}[1]
  {\renewcommand\theinnercustomthm{#1}\innercustomthm}
  {\endinnercustomthm}
\let\oldnl\nl
\newcommand{\nonl}{\renewcommand{\nl}{\let\nl\oldnl}}
\newcommand{\eproof}{\raisebox{-1mm}{\rule{0.5em}{1em}}}

\title{Make Workers Work Harder:
Decoupled Asynchronous Proximal Stochastic Gradient Descent}

\author{Yitan Li, Linli Xu, Xiaowei Zhong, Qing Ling \\
University of Science and Technology of China \\
etali@mail.ustc.edu.cn, linlixu@ustc.edu.cn, \{xwzhong,qingling\}@mail.ustc.edu.cn}

\begin{document}

\maketitle

\begin{abstract}
  Asynchronous parallel optimization algorithms for solving large-scale machine
  learning problems have drawn significant attention from academia to industry
  recently. This paper proposes a novel algorithm, decoupled asynchronous
  proximal stochastic gradient descent (DAP-SGD), to minimize an
  objective function that is the composite of the average of multiple empirical
  losses and a regularization term. Unlike the traditional asynchronous
  proximal stochastic gradient descent (TAP-SGD) in which the
  master carries much of the computation load, the proposed algorithm off-loads the
  majority of computation tasks from the master to workers, and leaves
  the master to conduct simple addition operations. This strategy
  yields an easy-to-parallelize algorithm, whose performance is
  justified by theoretical convergence analyses. To be specific,
  DAP-SGD achieves an $O(\log T/T)$ rate when the step-size is
  diminishing and an ergodic $O(1/\sqrt{T})$ rate when the step-size is
  constant, where $T$ is the number of total iterations.
\end{abstract}

\section{Introduction}

A majority of classical machine learning tasks can be formulated
as solving a general regularized optimization problem:
\begin{equation}\label{eqn:obj}
\begin{split}
\min_{\mathbf{x}\in \mathbb{R}^m} P(&\mathbf{x}) = f(\mathbf{x}) + h(\mathbf{x})\ , \\
\text{where}\quad &f(\mathbf{x}) \triangleq
\frac{1}{n}\sum_{i=1}^n f_i(\mathbf{x})\ .
\end{split}
\end{equation}
Given $n$ samples, $f_i(\mathbf{x})$ represents the empirical loss
of the $i^{th}$ sample with regard to the decision variable
$\mathbf{x}$, and $h(\mathbf{x})$ corresponds to a (usually
non-smooth) regularization term. Our goal is to find the optimal
solution, defined as $\mathbf{x}^*$, which minimizes the summation
of the averaged empirical loss and the regularization term over
the whole dataset.

With the enormous growth of data size $n$ and model complexity,
asynchronous parallel
algorithms~\cite{niu2011hogwild,duchi2011ddsc,mu2014ps,lian2015dl,ruiliang2014admm,FeyzmahdavianAJ14}
have become an important tool and received significant successes
for solving large scale machine learning problems in the form of
\eqref{eqn:obj}. Asynchronous parallel algorithms distribute
computation on multi-core systems~(shared memory architecture) or
multi-machine system~(parameter server architecture), whose
computation power generally scales up with the increasing number
of cores or machines. As a consequence, effective design and
implementation of asynchronous parallel algorithms is
critical for large scale machine learning. 

Numerous efforts have been devoted to this topic. Among them,
asynchronous stochastic gradient descent is proposed
in~\cite{niu2011hogwild,duchi2011ddsc}, and its performance is
guaranteed by theoretical convergence analyses. An asynchronous
proximal gradient descent algorithm is designed on the parameter
server architecture in~\cite{mu2014ps} with a distributed
optimization software provided. Convergence rate of asynchronous
stochastic gradient descent with a non-convex objective is
analyzed in \cite{lian2015dl}. Apart from work on asynchronous
gradient descent and its proximal variant, much attention has also
been attracted to asynchronous alternating direction method of
multipliers~(ADMM)~\cite{ruiliang2014admm}, asynchronous
stochastic coordinate
ascent~\cite{Ji2015sca,Ji2015sca2,olivier2015pcd,marevcek2015distributed,hong2014distributed,ZhouYDLX16}
and asynchronous dual stochastic coordinate
ascent~\cite{hsieh2015asdca}.

The traditional asynchronous proximal stochastic gradient method
(TAP-SGD) that solves \eqref{eqn:obj} works as follows. The
workers (multiple cores or machines) access samples, compute the
gradients of their corresponding empirical losses, and send to the
master. The master fuses the gradients and runs a proximal step on
the regularization term (more details are given in Section 2).
However, the performance of this paradigm is restricted when the
proximal operator is not an element-wise operation. For this case,
running proximal steps can be time-consuming, and the computation
in the master becomes the bottleneck of the whole system. We note
that this is common for many popular regularization terms, as
shown in Section 2. To avoid this difficulty, one has to design a
customized parallel computation for every single regularization
term, which makes the framework inflexible. For the
sake of speeding up computation and simplifying algorithm design,
we expect to design an alternative algorithm that is easier to
parallelize.

In light of this issue, this paper develops a decoupled
asynchronous proximal stochastic gradient descent (DAP-SGD), which
off-loads the majority of computation tasks (especially the
proximal steps) from the master to workers, and leaves the master
to conduct simple addition operations. This algorithmic framework
is suitable for many master/worker architectures including the
single machine multi-core system~(shared memory architecture)
where the master is the parameter updating thread and the workers
correspond to other threads processing samples, and the
multi-machine system~(parameter server architecture) where the
master is the central machine for storing and updating parameters
and the workers represent those machines for storing and
processing samples.

The main contributions of this paper are
highlighted as follows:
\begin{itemize}
    \item The proposed DAP-SGD algorithm off-loads the computation bottleneck from the master to workers. To be more specific,
    DAP-SGD allows workers to evaluate the proximal operators~(work harder) and the master only needs to do element-wise
    addition operations, which is easy to parallelize.
    \item Convergence analysis is provided for DAP-SGD. DAP-SGD achieves an $O(\log T/T)$ rate when the step-size is
    diminishing and an ergodic $O(1/\sqrt{T})$ rate when the step-size is constant, where $T$ is the number of total iterations.
\end{itemize}

\section{Traditional Asynchronous Proximal Stochastic Gradient Descent (TAP-SGD)}

We start from the synchronous proximal stochastic gradient descent
(P-SGD) algorithm that solves \eqref{eqn:obj}. P-SGD only requires
the gradient of one sample in a single iteration. Hence in large
scale optimization problems, it is a preferred surrogate for
proximal gradient descent \cite{beck2009fast,parikh2014proximal},
which requires computing gradients of all samples in a single
iteration. The recursion of P-SGD is
\begin{equation}\label{eqn:psgd}
\begin{split}
\mathbf{x}_{t+1} &= \text{Prox}_{\eta_t,h} (\mathbf{x}_t -
\eta_t~\triangledown\!f_{i_t}(\mathbf{x}_t))\ ,
\end{split}
\end{equation}
where $\text{Prox}_{\eta,h}(\mathbf{x}) = \arg\!\min_{\mathbf{y}}
\|\mathbf{y} - \mathbf{x}\|_2^2 /(2\eta) + h(\mathbf{y})$ denotes
a proximal operator, while $\eta_t$ is the step-size and $i_t$ is
the index of the selected sample in the $t^{th}$ iteration.

The traditional asynchronous proximal stochastic gradient descent
(TAP-SGD) algorithm is an asynchronous variant of P-SGD, as
summarized in Algorithm~\ref{alg:apsgd}. The master is the main
updating processor, while the workers provide the gradients of the
samples. Every worker receives the parameter (namely, decision
variables) $\mathbf{x}$ from the master, computes the gradient of
one random sample $\triangledown\!f_i (\mathbf{x})$ and sends it
to the master. Obviously, when one worker is computing and sending
its gradient, the master may update the parameter using the
gradients sent by the other workers in the previous time period.
As a consequence, the gradients received at the master are often
delayed, causing the main difference between P-SGD and TAP-SGD. In
the master, the delayed gradient received at the $t^{th}$
iteration is denoted by
$\triangledown\!f_{i_t}(\mathbf{x}_{d(t)})$ where $i_t$ indexes
the selected sample, $\mathbf{x}_{d(t)}$ refers to that the
parameter is the one from the $d(t)^{th}$ iteration, and $d(t)\in
[t-\tau, t]$ where $\tau$ stands for the maximum delay of the
system. Therefore, we can write the recursion of TAP-SGD as
\begin{equation}\label{eqn:apsgd}
\begin{split}
\mathbf{x}_{t+1} &= \text{Prox}_{\eta_t,h} (\mathbf{x}_t - \eta_t
\triangledown\!f_{i_t}(\mathbf{x}_{d(t)}))\ .
\end{split}
\end{equation}

\begin{algorithm}[ht]
\caption{Asynchronous Proximal Stochastic Gradient
Descent~(AP-SGD)} \label{alg:apsgd} \KwIn{Initialization
$\mathbf{x}_0$, $t = 0$, dataset with $n$ samples in which the
loss function of the $i^{th}$ sample is denoted by
$f_i(\mathbf{x})$, regularization term $h(\mathbf{x})$, maximum
number of iterations $T$, number of workers $S$, step-size in the
$t^{th}$ iteration $\eta_t$, maximum delay $\tau$}
\KwOut{$\mathbf{x}_T$}
  \nonl \textbf{Procedure of each worker $s\in [1,...,S]$} \DontPrintSemicolon\;
  \Repeat{procedure of master ends}{\PrintSemicolon
    Uniformly sample $i$ from $[1,...,n]$\;
    Obtain the parameter $\mathbf{x}$ from the master~(shared memory or parameter server)\;
    Evaluate the gradient of the $i^{th}$ sample over parameter $\mathbf{x}$, denoted by $\triangledown\!f_i(\mathbf{x})$\;
    Send $\triangledown\!f_i(\mathbf{x})$ to the master\;
  }
  \nonl \textbf{Procedure of master} \DontPrintSemicolon\;
  \setcounter{AlgoLine}{0}\For{$t = 0$ to $T-1$}{\PrintSemicolon
    Get a gradient $\triangledown\!f_{i_t}(\mathbf{x}_{d(t)})$~(the delay $t - d(t)$ is bounded by $\tau$)\;
    Update the parameter with the proximal operator
    $\mathbf{x}_{t+1} = \text{Prox}_{\eta_t,h} (\mathbf{x}_t - \eta_t \triangledown\!f_{i_t}(\mathbf{x}_{d(t)}))$\;
    $t = t + 1$\;
  }
\end{algorithm}

Observe that the updating procedure of the master is the
computational bottleneck of the TAP-SGD algorithm. When the
proximal step is time-consuming to calculate, the workers must
wait for a long time to receive updated parameters, which
significantly degrades the performance of the system. To avoid
this difficulty, one has to design a customized parallel
computation for every single regularization term, which makes
the framework inflexible. In a multi-machine system with multiple
masters, such parallelized proximal operators will also cause
complicated network communications between masters.

\subsection*{Coupled Proximal Operators}

In practice, many widely used (usually non-smooth) regularization
terms are associated with coupled proximal operators, which lead
to high computational complexity, including group lasso
regularization~\cite{friedman2010note}, fused lasso
regularization~\cite{Jun2010fusedlasso}, nuclear norm
regularization~\cite{Shuiwang2009nuclear,Cai2010svt}, etc.

\textbf{The proximal operator of group lasso regularization} $h(\mathbf{x}) = \lambda \sum_{i=1}^{g}\|\mathbf{x}_{k_i:(k_{i+1}-1)}\|_2$:
\begin{equation}\label{eqn:prox_grouplasso}
\begin{split}
\text{Prox}_{\eta,h}(\mathbf{x}) = \arg\!\min_{\mathbf{y}}& \frac{1}{2\eta} \|\mathbf{y} - \mathbf{x}\|_2^2 + \lambda \sum_{i=1}^{g}\|\mathbf{y}_{k_i:(k_{i+1}-1)}\|_2\ . \\
\end{split}
\end{equation}
Here $g$ is the number of groups and $k_1 = 1 < ...\ k_i < k_{i+1}
...< k_{g+1} = m+1$. The closed-form solution of the proximal
operator above is
\begin{equation}
\begin{split}
{[\text{Prox}_{\eta,h}(\mathbf{x})]}_{k_i:(k_{i+1}-1)} = \mathbf{x}_{k_i:(k_{i+1}-1)} \left(1 - \frac{\lambda}{\|\mathbf{x}_{k_i:(k_{i+1}-1)}\|_2} \right)_+ \ . \\
\end{split}
\end{equation}

For the group lasso regularization, the proximal
operator~(\ref{eqn:prox_grouplasso}) is separated into $g$ groups.
When partitions of groups are unbalanced, it will be hard to speed
up the computation with parallelization.

\textbf{The proximal operator of simplified fused lasso regularization} $h(\mathbf{x}) = \lambda \sum_{i=1}^{m-1}\|\mathbf{x}_i-\mathbf{x}_{i+1}\|_1$:
\begin{equation}\label{eqn:prox_fusedlasso}
\begin{split}
\text{Prox}_{\eta,h}(\mathbf{x}) &= \arg\!\min_{\mathbf{y}} \frac{1}{2\eta} \|\mathbf{y} - \mathbf{x}\|_2^2 + \lambda \sum_{i=1}^{m-1}\|\mathbf{y}_i-\mathbf{y}_{i+1}\|_1 \\
&= \mathbf{y} - \mathbf{R}^T \mathbf{z}^* \ , 
\end{split}
\end{equation}
\begin{equation*}
\begin{split}
\text{where} \quad &\mathbf{R} =
\begin{bmatrix} 1 & -1 & 0 & ... & 0 \\
0 & 1 & -1 & ... & 0 \\
&& ... && \\
0 & ... & 0 & 1 & -1 \end{bmatrix} \in \mathbb{R}^{(m-1)\times m}, \\
& \mathbf{z}^* = \arg\!\min_{\|\mathbf{z}\|_{\infty}\leq \eta\lambda} \frac{1}{2}\|\mathbf{R}^{\mathrm{T}}\mathbf{z}\|_2^2 - < \mathbf{R}^{\mathrm{T}}\mathbf{z},\mathbf{y} > \ . \\
\end{split}
\end{equation*}
For the simplified fused lasso regularization, the proximal
operator~(\ref{eqn:prox_fusedlasso}) has a closed form solution.
However, solving $\mathbf{z}^*$ involves a subproblem that is
time-consuming.

\textbf{The proximal operator of nuclear norm regularization} $h(\mathbf{X}) = \lambda \|\mathbf{X}\|_*$:
\begin{equation}\label{eqn:prox_nuclear}
\begin{split}
\text{Prox}_{\eta,h}(\mathbf{X}) &= \arg\!\min_{\mathbf{Y}} \frac{1}{2\eta} \|\mathbf{Y} - \mathbf{X}\|_F^2 + \lambda \|\mathbf{Y}\|_* \\
&= \mathbf{U} \mathbf{\hat{\Sigma}} \mathbf{V}^{\mathrm{T}}\ ,
\end{split}
\end{equation}
where $\mathbf{X} = \mathbf{U}\mathbf{\Sigma}
\mathbf{V}^{\mathrm{T}}$ calculated from singular value
decomposition, $\sigma_i$ is the $i^{th}$ singular value of
$\mathbf{X}$, $\hat{\sigma}_i=\max(\sigma_i-\eta\lambda, 0)$ is
the $i^{th}$ element of $\hat{\bm{\sigma}})$, and
$\hat{\mathbf{\Sigma}}= \mathrm{Diag}(\hat{\bm{\sigma}})$. For the
nuclear norm regularization, the proximal
operator~(\ref{eqn:prox_nuclear}) involves singular value
decomposition, which is challenging especially for large scale
problems.

As discussed above, evaluating the proximal operator can be a
computational bottleneck and limits the performance of TAP-SGD.
This motivates us to design a novel asynchronous parallel
algorithm, which decouples and distributes the calculation of the
proximal operator to the workers.

\section{Decoupled Asynchronous Proximal Stochastic Gradient Descent (DAP-SGD)}

The key idea of the decoupled asynchronous proximal stochastic
gradient descent~(DAP-SGD) algorithm is to off-load the
computational bottleneck from the master to the workers. The
master no longer takes care of the proximal operators; instead, it
only needs to conduct element-wise addition operations. On the
other hand, the workers must work harder: they evaluate the
proximal operators independently, without caring about the
parallel mechanism.

The procedure of DAP-SGD is summarized in Algorithm
\ref{alg:dapsgd}. Each worker evaluates the proximal operator and
sends update information (namely, innovation) $\Delta =
\mathbf{x}' - \mathbf{x}$ to the master. In the master, the
delayed update information $\Delta_{d(t)} = \mathbf{x}'_{d(t)} -
\mathbf{x}_{d(t)}$ is used to modify the parameter $\mathbf{x}$.
Obviously, parameter updating in the master is no longer the
computational bottleneck of the system, since it only involves
element-wise addition operations.

\begin{algorithm}[!htb]
\caption{Decoupled Asynchronous Proximal Stochastic Gradient
Descent~(DAP-SGD)} \label{alg:dapsgd} \KwIn{Initialization
$\mathbf{x}_0$, $t = 0$, dataset with $n$ samples in which loss
function of the $i^{th}$ sample is denoted by $f_i(\mathbf{x})$,
regularization term $h(\mathbf{x})$, maximum number of iterations
$T$, number of workers $S$, step-size in the $t^{th}$ iteration
$\eta_t$, maximum delay $\tau$} \KwOut{$\mathbf{x}_T$}
  \nonl \textbf{Procedure of each worker $s\in [1,...,S]$} \DontPrintSemicolon\;
  \Repeat{procedure of master end}{\PrintSemicolon
    Uniformly sample $i$ from $[1,...,n]$\;
    Obtain parameter $\mathbf{x}$ and step-size $\eta$ from master~(shared memory or parameter server)\;
    Evaluate the gradient of the $i^{th}$ sample over parameter $\mathbf{x}$, denoted by $\triangledown\!f_i(\mathbf{x})$\;
    Evaluate the proximal operator $\mathbf{x}' = \text{Prox}_{\eta,h} (\mathbf{x} - \eta ~\triangledown\!f_i(\mathbf{x}))$\;
    Send update information $\Delta = \mathbf{x}' - \mathbf{x}$ to the master\;
  }
  \nonl \textbf{Procedure of master} \DontPrintSemicolon\;
  \setcounter{AlgoLine}{0}\For{$t = 0$ to $T-1$}{\PrintSemicolon
    Get $\Delta_{d(t)} = \mathbf{x}'_{d(t)} - \mathbf{x}_{d(t)}$ from one worker~(the delay $t - d(t)$ is bounded by $\tau$)\;
    Update parameter with
    $\mathbf{x}_{t+1} = \mathbf{x}_t + \Delta_{d(t)}$\;
    $t = t + 1$\;
  }
\end{algorithm}

The recursion of DAP-SGD is
\begin{equation}\label{eqn:dapsgd}
\begin{split}
\mathbf{x}'_{d(t)} &= \text{Prox}_{\eta,h} (\mathbf{x}_{d(t)} - \eta_{d(t)} \triangledown\!f_{i_{d(t)}}(\mathbf{x}_{d(t)})) \ , \\
\mathbf{x}_{t+1} &= \mathbf{x}_t + \mathbf{x}'_{d(t)} -
\mathbf{x}_{d(t)} \ .
\end{split}
\end{equation}
Comparing the recursions of TAP-SGD~(\ref{eqn:apsgd}) and
DAP-SGD~(\ref{eqn:dapsgd}), we can observe that the DAP-SGD
recursion~(\ref{eqn:dapsgd}) splits the proximal operator and
parameter updating step\footnote{Note that both TAP-SGD and
DAP-SGD can support mini-batch updating.}. This is the why we call
the proposed algorithm ``\textbf{decoupled}''. The benefit of
decoupling is that the computational bottleneck (for example, the
unbalanced partitioned groups in (\ref{eqn:prox_grouplasso}), the
subproblem in (\ref{eqn:prox_fusedlasso}), and the
singular value decomposition in (\ref{eqn:prox_nuclear})) no
longer lies in the master. The workers conduct these operations,
which improves the performance of the system. Below, we further
analyze the convergence properties of DAP-SGD theoretically.

\section{Convergence Analysis}

This section gives theorems that establish the convergence
properties of DAP-SGD. The detailed proofs are presented in the
appendix. We start from some basic assumptions.

The first two assumptions are about the properties of the averaged
empirical cost $f(\mathbf{x})$.

\begin{assumption}\label{asp:lcg}
\textbf{Lipschitz continuous gradient of
$\triangledown\!f(\mathbf{x})$}: The function $f(\mathbf{x})$ is
differentiable and its gradient $\triangledown\!f(\mathbf{x})$ is
Lipschitz continuous with constant $L$. Namely, the following two
equivalent inequalities hold:
\begin{equation}\label{eqn:lcg}
\begin{split}
\ f(\mathbf{x})\leq f(\mathbf{y}) + \left<
\triangledown\!f(\mathbf{y}),\mathbf{x} - \mathbf{y} \right> +
\frac{L}{2}\|\mathbf{x} - \mathbf{y}\|_2^2\ , \quad \forall\
\mathbf{x},\mathbf{y},
\end{split}
\end{equation}
and
\begin{equation}\label{eqn:lcg2}
\begin{split}
\ \frac{1}{L}\|\triangledown\!f(\mathbf{x}) -
\triangledown\!f(\mathbf{y})\|^2\leq \left<
\triangledown\!f(\mathbf{x}) -
\triangledown\!f(\mathbf{y}),\mathbf{x}-\mathbf{y} \right> \leq
L\|\mathbf{x}-\mathbf{y}\|^2\ , \quad \forall\
\mathbf{x},\mathbf{y}.
\end{split}
\end{equation}
\end{assumption}

\begin{assumption}\label{asp:sc}
\textbf{Strong convexity of $f(\mathbf{x})$}: The function
$f(\mathbf{x})$ is strongly convex with constant $\mu$. Namely,
the following inequality holds:
\begin{equation}\label{eqn:sc}
\begin{split}
\ f(\mathbf{x}) \geq f(\mathbf{y}) + \left<
\triangledown\!f(\mathbf{y}),\mathbf{x} - \mathbf{y} \right> +
\frac{\mu}{2}\|\mathbf{x} - \mathbf{y}\|_2^2 \ , \quad \forall\
\mathbf{x},\mathbf{y}.
\end{split}
\end{equation}
\end{assumption}

The next assumption bounds the variance of sampling a random
gradient $\triangledown\!f_i(\mathbf{x})$ to replace the true
gradient $\triangledown\!f(\mathbf{x})$.

\begin{assumption}\label{asp:varb}
\textbf{Bounded variance of gradient evaluation}: The variance of
a selected gradient is bounded by a constant $C_f$:
\begin{equation}\label{eqn:varb}
\begin{split}
\ \mathbb{E} \|\triangledown\!f_i(\mathbf{x}) -
\triangledown\!f(\mathbf{x})\|_2^2 \leq C_f \ , \quad \forall\
\mathbf{x}.
\end{split}
\end{equation}
\end{assumption}

The last two assumptions are about the properties of the
regularization term $h(\mathbf{x})$.

\begin{assumption}\label{asp:sc-h}
\textbf{Convexity of $h(\mathbf{x})$}: The function
$h(\mathbf{x})$ is convex. Namely, the following inequality holds:
\begin{equation}\label{eqn:sc-h}
\begin{split}
\ h(\mathbf{x}) \geq h(\mathbf{y}) + \left<
\partial h(\mathbf{y}),\mathbf{x} - \mathbf{y} \right>\ , \quad \forall\
\mathbf{x},\mathbf{y},
\end{split}
\end{equation}
\end{assumption}
where $\partial h(\mathbf{x})$ stands for any subgradient of
$h(\mathbf{x})$.

\begin{assumption}\label{asp:subb}
\textbf{Bounded subgradient of $h(\mathbf{x})$}: The squared
subgradient of $h(\mathbf{x})$ is bounded by a constant $C_h$
\begin{equation}\label{eqn:subb}
\begin{split}
\|\partial h(\mathbf{x})\|_2^2 \leq C_h \ .
\end{split}
\end{equation}
\end{assumption}

An immediate result from Assumption \ref{asp:subb} is that,
$\triangledown\!f(\mathbf{x^*})$ is also bounded where $x^*$ is
the optimal solution to \eqref{eqn:obj}, as given in the following
corollary.

\begin{corollary}\label{asp:optb}
\textbf{Bounded gradient of $f(\mathbf{x})$ at the optimum}: Let
$\mathbf{x}^* = \arg\!\min_{\mathbf{x}} f(\mathbf{x}) +
h(\mathbf{x})$ be the optimal solution to \eqref{eqn:obj}, then we
have
\begin{equation}\label{eqn:optb}
\begin{split}
\|\triangledown\!f(\mathbf{x}^*)\|_2^2 = \|\partial
h(\mathbf{x}^*)\|_2^2 \leq C_h \ .
\end{split}
\end{equation}
\end{corollary}

Assumptions \ref{asp:lcg}, \ref{asp:sc}, \ref{asp:varb} and
\ref{asp:sc-h} are common in the convergence analysis of
stochastic gradient descent
algorithms~\cite{niu2011hogwild,mu2014ps,lian2015dl,nesterov2013introductory,guanghui2009sgd}.
Assumption \ref{asp:subb} is due to the (usually non-smooth)
regularization term $h(\mathbf{x})$, and is reasonable for many
non-smooth regularization terms such as $L_1$ regularization,
group lasso, fused lasso and nuclear norm, etc. Next we provide
the constant upper bounds of subgradients for these non-smooth
regularization terms. In the following part, $\partial$ denotes
the set of subderivatives, and with a slight abuse of notation,
also denotes any element (namely, subgradient) in the set.

\textbf{Upper bound of subgradient for $L_1$ regularization}
$\|x\|_1$:
\begin{equation}\label{eqn:ub_l1}
\begin{split}
\|\partial\|\mathbf{x}\|_1\|_2 \leq m \ .
\end{split}
\end{equation}

\textbf{Upper bound of subgradient for group lasso regularization}
$\sum_{i=1}^{g}\|\mathbf{x}_{k_i:(k_{i+1}-1)}\|_2$:
\begin{equation}\label{eqn:ub_grouplasso}
\begin{split}
&\quad\quad\quad\quad\quad \ \left\|\partial \sum_{i=1}^{g}\|\mathbf{x}_{k_i:(k_{i+1}-1)}\|_2\right\|\leq g\ , \\
&\text{where} \quad \ \partial\|\mathbf{x}_{k_i:(k_{i+1}-1)}\|_2 =
\begin{cases}
\frac{1}{\|\mathbf{x}_{k_i:(k_{i+1}-1)}\|} \mathbf{x}_{k_i:(k_{i+1}-1)} & \text{if}\ \mathbf{x}_{k_i:(k_{i+1}-1)}\neq 0, \\
\{g|\|g\|_2\leq 1\} & \text{if}\ \mathbf{x}_{k_i:(k_{i+1}-1)} = 0 \ .
\end{cases}
\end{split}
\end{equation}

\textbf{Upper bound of subgradient for simplified fused lasso
regularization}
$\sum_{i=1}^{m-1}\|\mathbf{x}_i-\mathbf{x}_{i+1}\|_1 =
\|\mathbf{R}\mathbf{x}\|_1$:
\begin{equation}\label{eqn:ub_fusedlasso}
\begin{split}
\|\partial \|\mathbf{R}\mathbf{x}\|_2\|_2 = & \|\mathbf{R}^\mathrm{T} SGN(\mathbf{R}\mathbf{x})\|_2 \leq \sum_i \|\mathbf{R}_{:,i}\|_2 \|SGN(\mathbf{R}\mathbf{x})\|_2\leq (m-1) \sum_i \|\mathbf{R}_{:,i}\|_2 \\
\leq & \sqrt{2}m(m-1),
\end{split}
\end{equation}
where $SGN$~\cite{Jun2010fusedlasso} is a function whose output is within $[-1,1]$.

\textbf{Upper bound of subgradient of nuclear norm regularization}
$\|\mathbf{X}\|_*, \mathbf{X}\in \mathbb{R}^{m\times q},d =
\min(m,q)$:
\begin{equation}\label{eqn:ub_nuclear}
\begin{split}
&\|\partial\|\mathbf{X}\|_*\|_F\leq
\|\mathbf{U}\mathbf{V}^\mathrm{T}\|_F + \|\mathbf{W}\|_F\leq
\|\mathbf{U}\|_F\|\mathbf{V}^\mathrm{T}\|_F + \|\mathbf{W}\|_F
\leq \mathrm{rank}(\mathbf{X})^2 + d \leq d^2 + d \ ,
\end{split}
\end{equation}
where $\partial\|\mathbf{X}\|_* =
\{\mathbf{U}\mathbf{V}^\mathrm{T} + \mathbf{W} | \mathbf{W}\in
\mathbb{R}^{m\times q}$, $\mathbf{U}^\mathrm{T}\mathbf{W} = 0$,
$\mathbf{W}\mathbf{V} = 0$, $\|\mathbf{W}\|_2\leq 1$, $\mathbf{X}
= \mathbf{U}\mathbf{\Sigma} \mathbf{V}^\mathrm{T}\}$.

Under the assumptions given above, we prove that DAP-SGD achieves
an $O(\log T/T)$ rate when the step-size is diminishing (Theorem
\ref{thm:sc}) and an ergodic $O(1/\sqrt{T})$ rate when the
step-size is constant (Theorem \ref{thm:scfixed}), where $T$ is
the number of total iterations. The proofs of the theorems are
given in the appendix.

\begin{theorem}\label{thm:sc}
Suppose that the cost function of (1) satisfies the following
conditions: $f(\mathbf{x})$ is strongly convex with constant $\mu$
and $h(\mathbf{x})$ is convex; $f(\mathbf{x})$ is differentiable
and $\triangledown\! f(\mathbf{x})$ is Lipschitz continuous with
constant $L$; $\mathbb{E} \|\triangledown\!f_i(\mathbf{x}) -
\triangledown\!f(\mathbf{x})\|_2^2 \leq C_f$; $\|\partial
h(\mathbf{x})\|_2^2 \leq C_h$. Define the optimal solution of (1)
as $\mathbf{x}^*$. At time $t$, set the step-size of the DAP-SGD
recursion (8) as $\eta_t = O(1/t)$. Then the iterate generated by
(\ref{eqn:dapsgd}) at time $T$, denoted by $\mathbf{x}_{T}$, satisfies
\begin{equation}
\begin{split}
\mathbb{E}\|\mathbf{x}_T - \mathbf{x}^*\|_2^2 \leq
O\left(\frac{\log T}{T}\right).
\end{split}
\end{equation}
\end{theorem}

\begin{theorem}\label{thm:scfixed}
Suppose that the cost function of (1) satisfies the following
conditions: $f(\mathbf{x})$ is strongly convex with constant $\mu$
and $h(\mathbf{x})$ is convex; $f(\mathbf{x})$ is differentiable
and $\triangledown\! f(\mathbf{x})$ is Lipschitz continuous with
constant $L$; $\mathbb{E} \|\triangledown\!f_i(\mathbf{x}) -
\triangledown\!f(\mathbf{x})\|_2^2 \leq C_f$; $\|\partial
h(\mathbf{x})\|_2^2 \leq C_h$. Define the optimal solution of (1)
as $\mathbf{x}^*$. At time $t$, fix the step-size of the DAP-SGD
recursion (8) $\eta_t$ as $\eta = O(1/\sqrt{T})$, where $T$ is the
maximum number of iterations. Define the iterate generated by (8)
at time $t$ as $\mathbf{x}_t$. Then the running average iterate
generated by (\ref{eqn:dapsgd}) at time $T$, denoted by
$\bar{\mathbf{x}}_T = \sum_{t = 0}^T \mathbf{x}_t/(T+1)$,
satisfies
\begin{equation}
\begin{split}
\mathbb{E}\|\bar{\mathbf{x}}_T - \mathbf{x}^*\|_2^2 \leq
O\left(\frac{1}{\sqrt{T}}\right).
\end{split}
\end{equation}
\end{theorem}

\section{Experiments}

\begin{figure}[t]\centering
\subfigure[$L_1$]{ \label{fig:exp1l1_time}
\includegraphics[width=0.23\columnwidth]{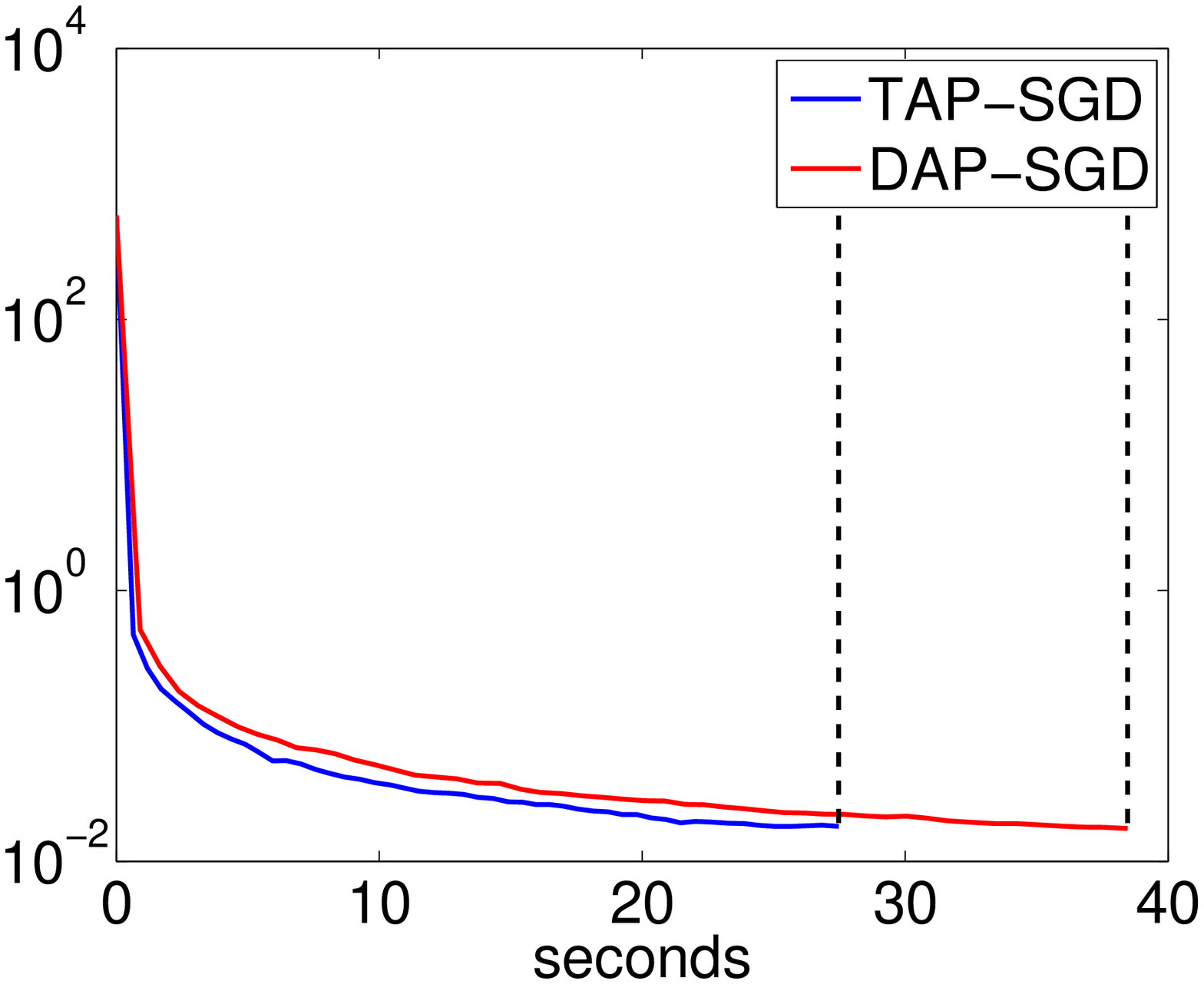}
} \subfigure[group lasso]{ \label{fig:exp1glasso_time}
\includegraphics[width=0.23\columnwidth]{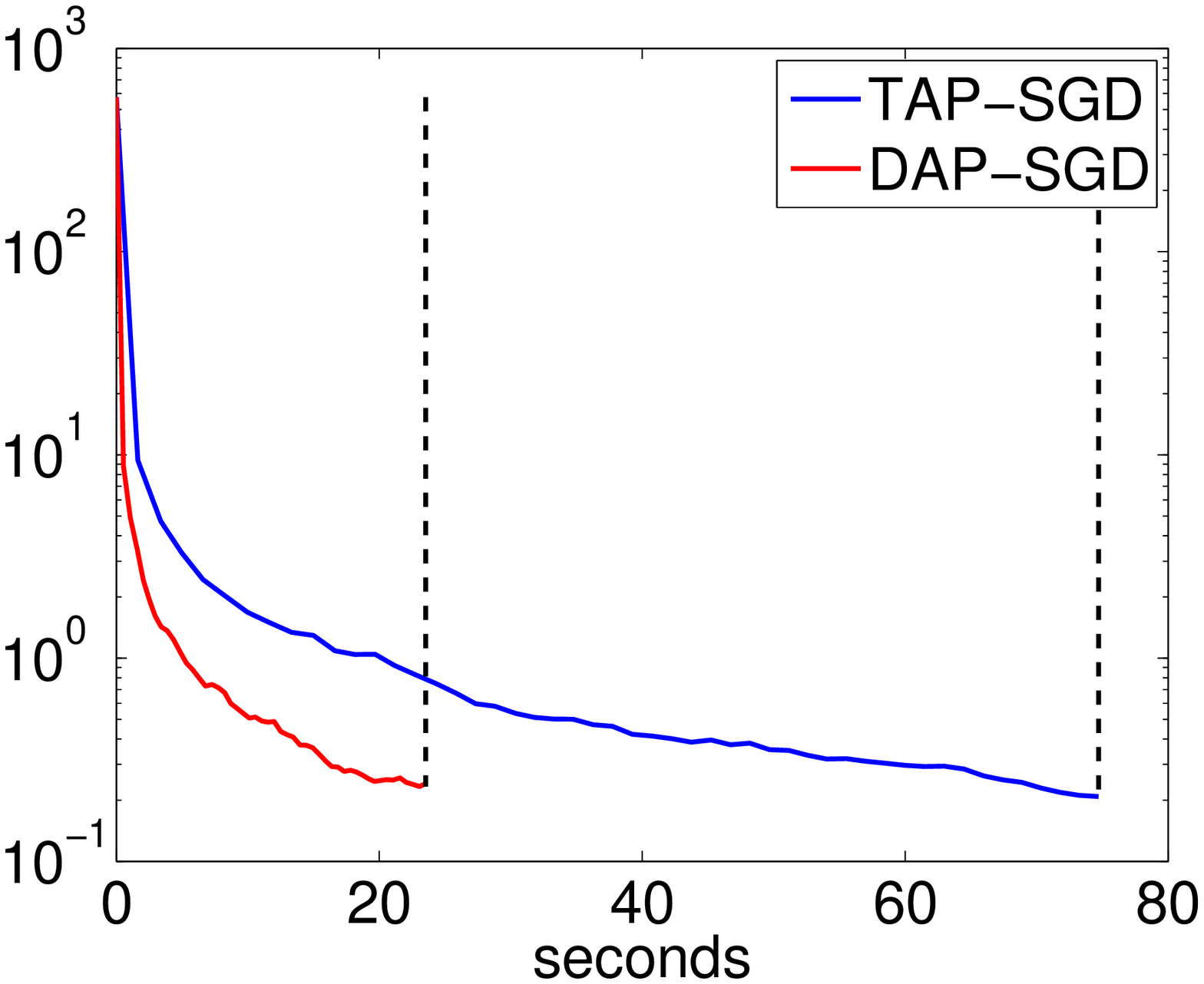}
} \subfigure[fused lasso]{ \label{fig:exp1sfusedlasso_time}
\includegraphics[width=0.23\columnwidth]{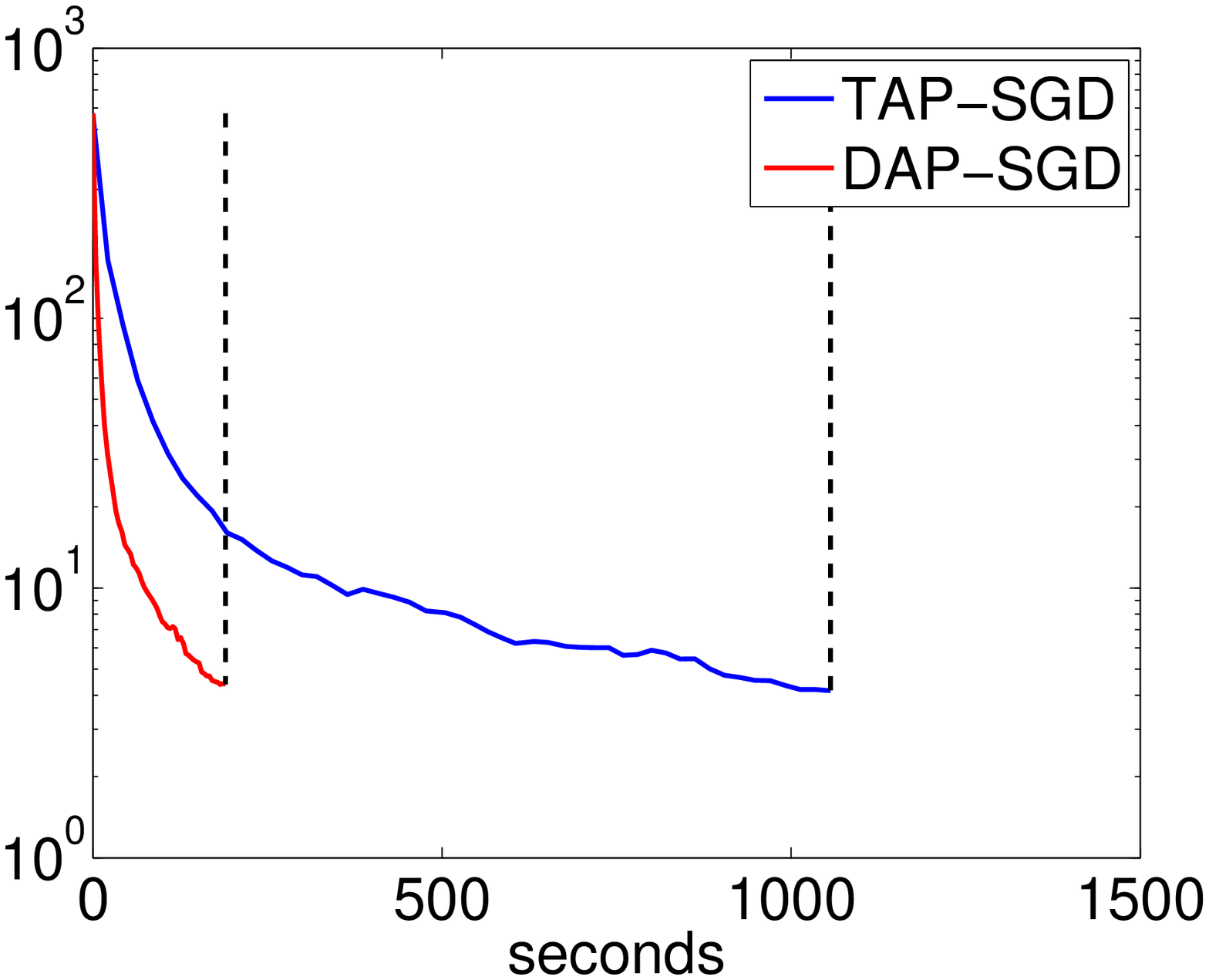}
} \subfigure[nuclear norm]{ \label{fig:exp1svt_time}
\includegraphics[width=0.23\columnwidth]{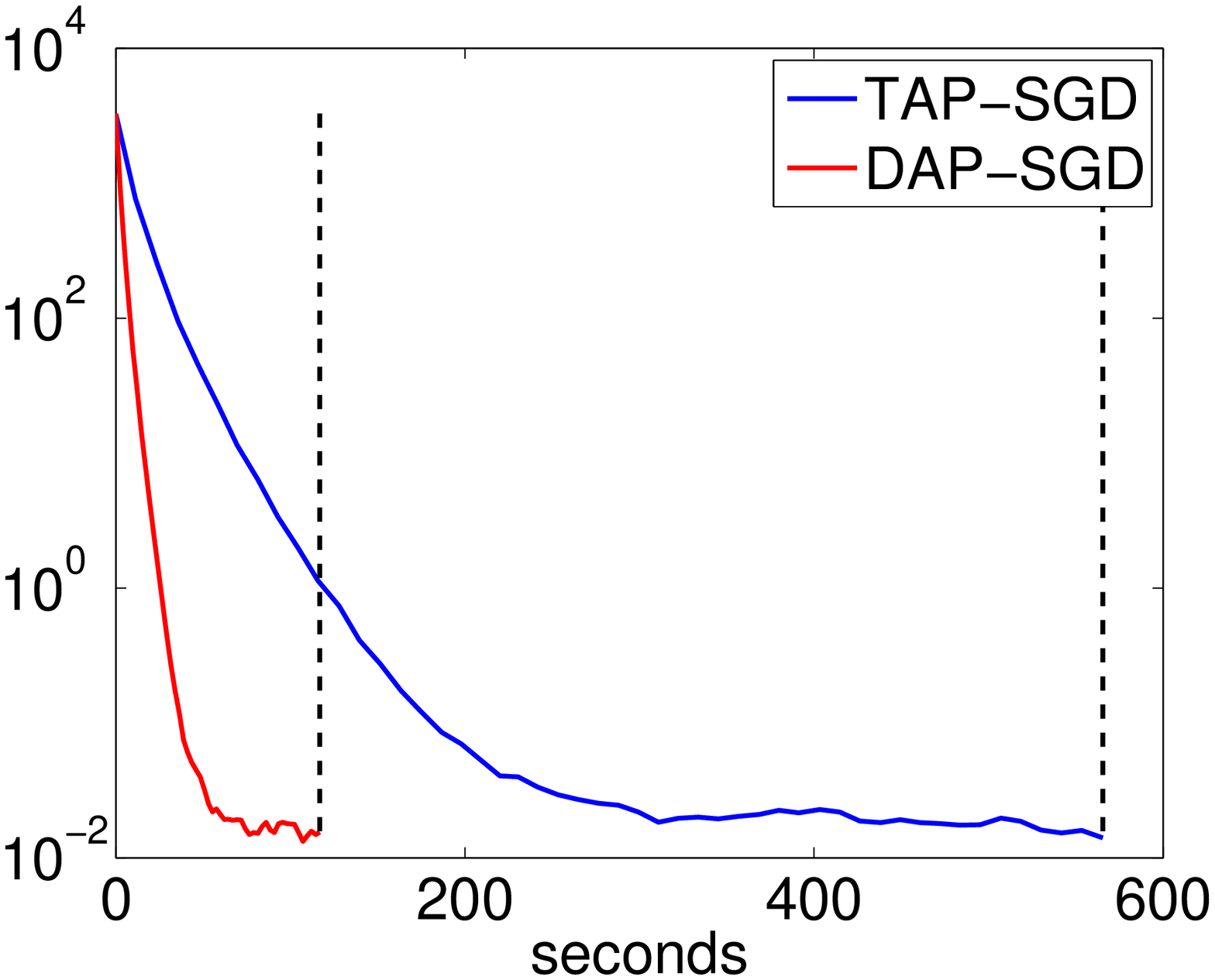}
}
\\
\subfigure[$L_1$]{ \label{fig:exp1l1_iter}
\includegraphics[width=0.23\columnwidth]{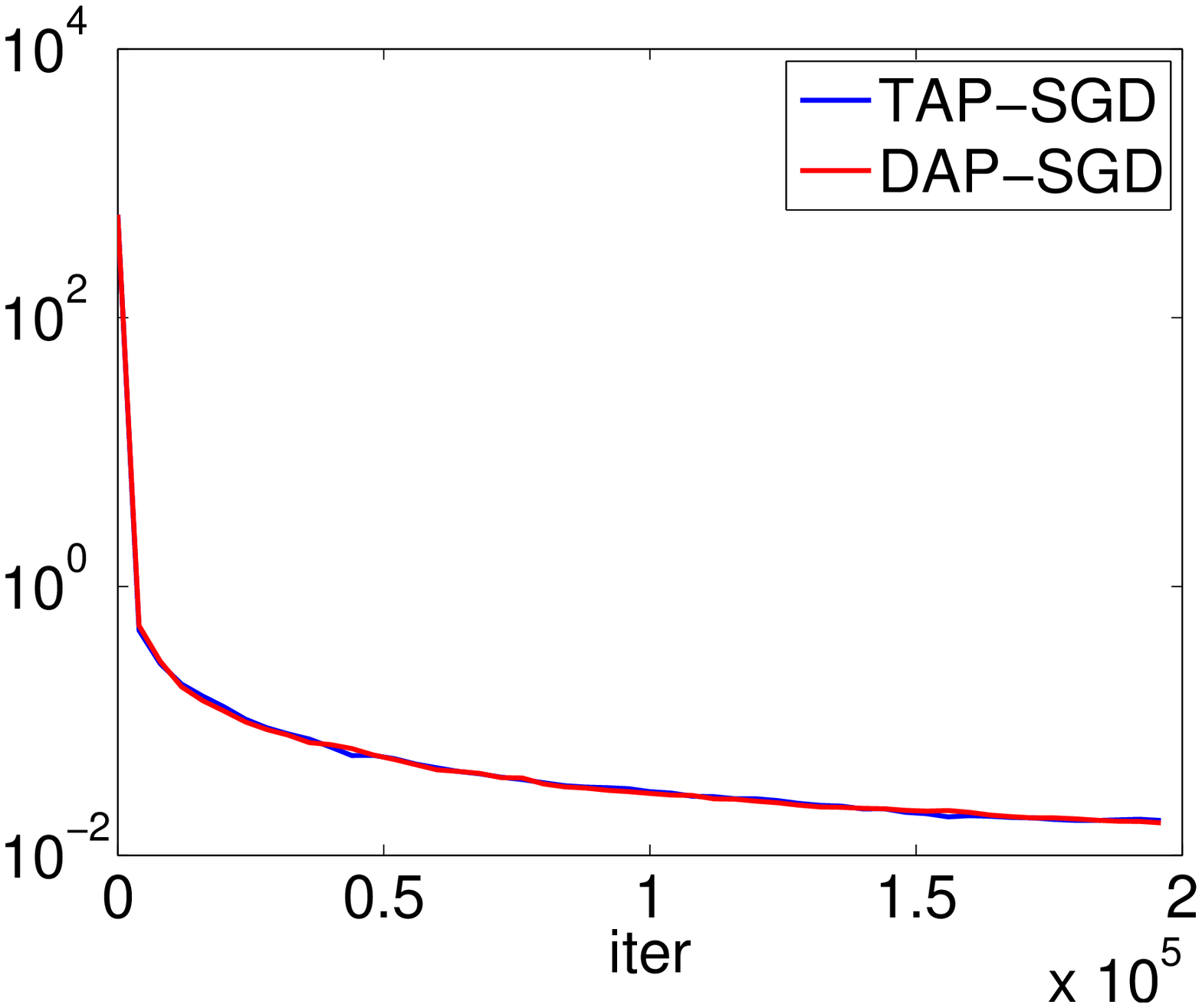}
} \subfigure[group lasso]{ \label{fig:exp1glasso_iter}
\includegraphics[width=0.23\columnwidth]{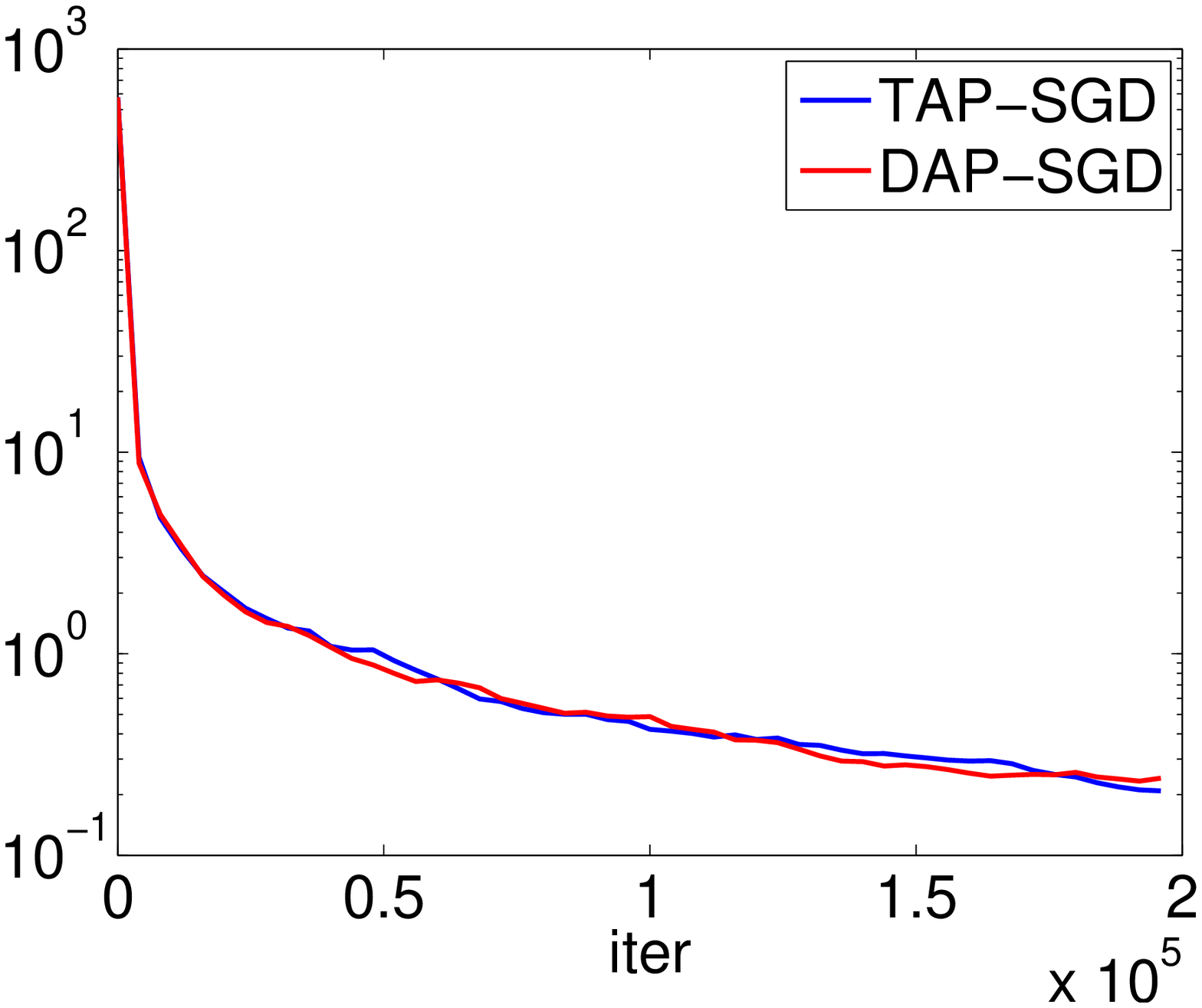}
} \subfigure[fused lasso]{ \label{fig:exp1sfusedlasso_time}
\includegraphics[width=0.23\columnwidth]{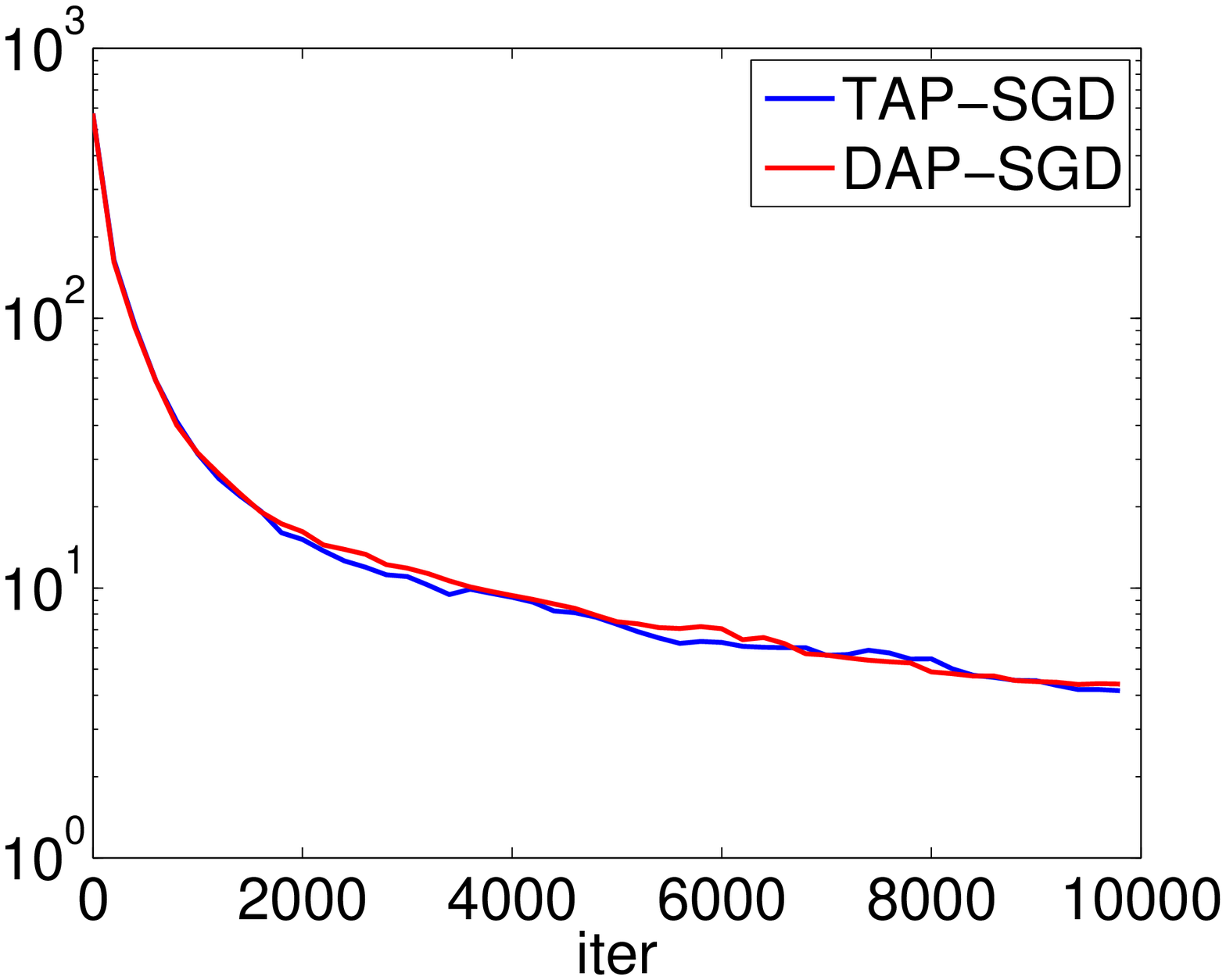}
} \subfigure[nuclear norm]{ \label{fig:exp1svt_time}
\includegraphics[width=0.23\columnwidth]{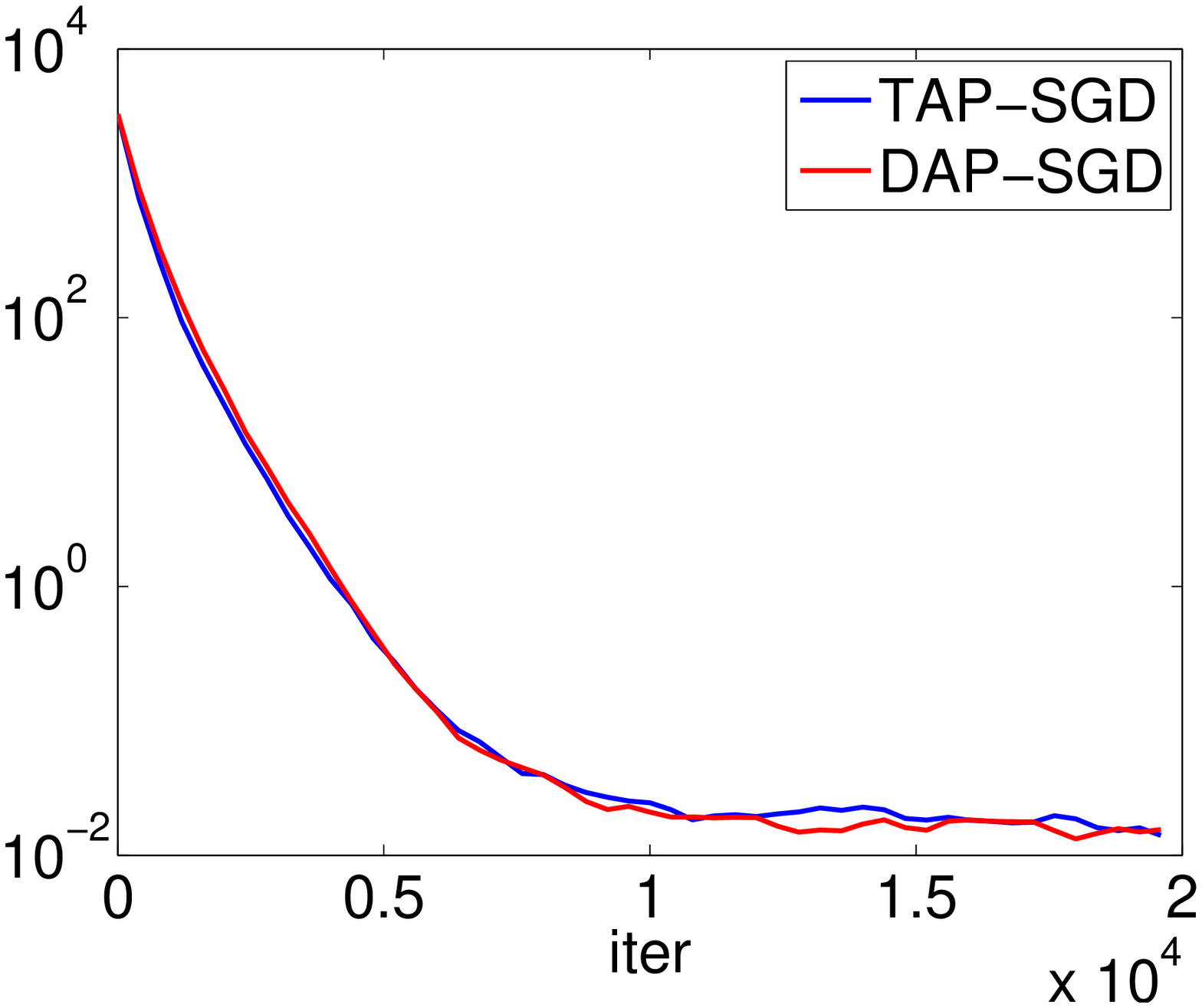}
}\caption{Comparison of TAP-SGD and DAP-SGD in terms of time and
number of iterations. The Y-axis shows the log distance between
the solution generated by an algorithm and the optimal solution,
denoted by $\log \|\mathbf{x} - \mathbf{x}^*\|_2^2$. Results of
$L_1$, group lasso, simplified fused lasso and nuclear norm
regularized objectives are shown in columns from left to right,
respectively. Top and bottom rows correspond to the results
regarding time and number of iterations,
respectively.}\label{fig:exp1}
\end{figure}

We compare the proposed DAP-SGD algorithm with TAP-SGD in a
consistent way without assuming the data is sparse. The
implementation is based on the single machine multi-core
system~(shared memory architecture). Both algorithms are
implemented in C++ and run on a multi-core server.
Singular value
decomposition~(SVD) is calculated by
eigen3\footnote{eigen.tuxfamily.org}. The parameters are locked
while they are being updated. The lock operation will slow down
the computation; however it guarantees that the implementation
conforms to the algorithm and its corresponding convergence
analysis.

Without loss of generality, we choose the least square loss with a
non-smooth regularization term as the optimization objective:
\begin{equation}
\begin{split}
\min_{\mathbf{x}\in \mathbb{R}^m} P(\mathbf{x}) = f(\mathbf{x}) +
h(\mathbf{x}) = \frac{1}{n}\sum_{i=1}^n
\left[\|\mathbf{x}^{\mathrm{T}}\mathbf{s}_i - y_i\|_2^2 + \lambda
\|\mathbf{x}\|_2^2\right] + h(\mathbf{x})
\end{split}.
\end{equation} In the case of nuclear norm regularization, the loss function
$f(\mathbf{x})$ becomes the multi-target least square loss
$f(\mathbf{X}) = \frac{1}{n}\sum_{i=1}^n
\left[\|\mathbf{X}^{\mathrm{T}}\mathbf{s}_i - \mathbf{y}_i\|_2^2 + \lambda \|\mathbf{X}\|_F^2 \right]$
correspondingly.

In the implementation TAP-SGD, the proximal operator of the $L_1$
regularized objective can be parallelized easily, while the
proximal operators of group lasso, simplified fused lasso and
nuclear norm are not parallelized due to their coupled and
non-element-wise operations. On the other hand, the procedure of
the master in the proposed DAP-SGD only involves simple
element-wise operations.

\textbf{Experimental Setup}. We conduct two experiments to
evaluate the algorithms with 4 different non-smooth regularization
terms~($L_1$, group lasso, simplified fused lasso, nuclear norm)
regarding the running time and number of iterations, as well as
the speedup. Data is generated randomly. In the first experiment,
for the 4 different objectives, the number of samples $n$ is set
to $1\times 10^3$, $1\times 10^3$, $1\times 10^3$, and $4\times
10^3$, while the length of the parameter is set to $5\times 10^3$,
$5\times 10^3$, $5\times 10^3$ and $2\times 10^3$ (in the form of
a $50\times 40$ matrix for nuclear norm regularization),
respectively. The number of iterations $T$ is set to $2\times
10^5$, $2\times 10^5$, $1\times 10^4$ and $2\times 10^4$, and the
step-size $\eta_t$ is set to $\frac{1}{2\times 10^5 + 200t}$,
$\frac{1}{2\times 10^5 + 200t}$, $\frac{1}{2\times 10^5 + 200t}$
and $\frac{1}{2\times 10^4 + t}$, respectively, which is
decreasing with iterations. The hyper-parameter $\lambda$ is set
to $200$, $200$, $200$, $0.1$ correspondingly. In the second
experiment of evaluating the speedup, the settings are identical
to the first experiment except that the number of iterations for
simplified fused norm and nuclear norm regularized objectives is
set to $10^4$ and $2\times 10^4$, and the number of parameters for
$L_1$ and group lasso regularized objectives is set to $5\times
10^4$. The total time cost of a system consists of two parts:
evaluation of updating information in the workers and updating in
the master. If we can speed up both with $k$ times, then we can
achieve a $k${-}speed up in the ideal case. In our experiment, the
number of updating threads running in parallel and maximum delay
$\tau$ in the master is fixed to the number of workers.

\begin{figure}[t]\centering
\subfigure[$L_1$]{ \label{fig:exp2l1}
\includegraphics[width=0.23\columnwidth]{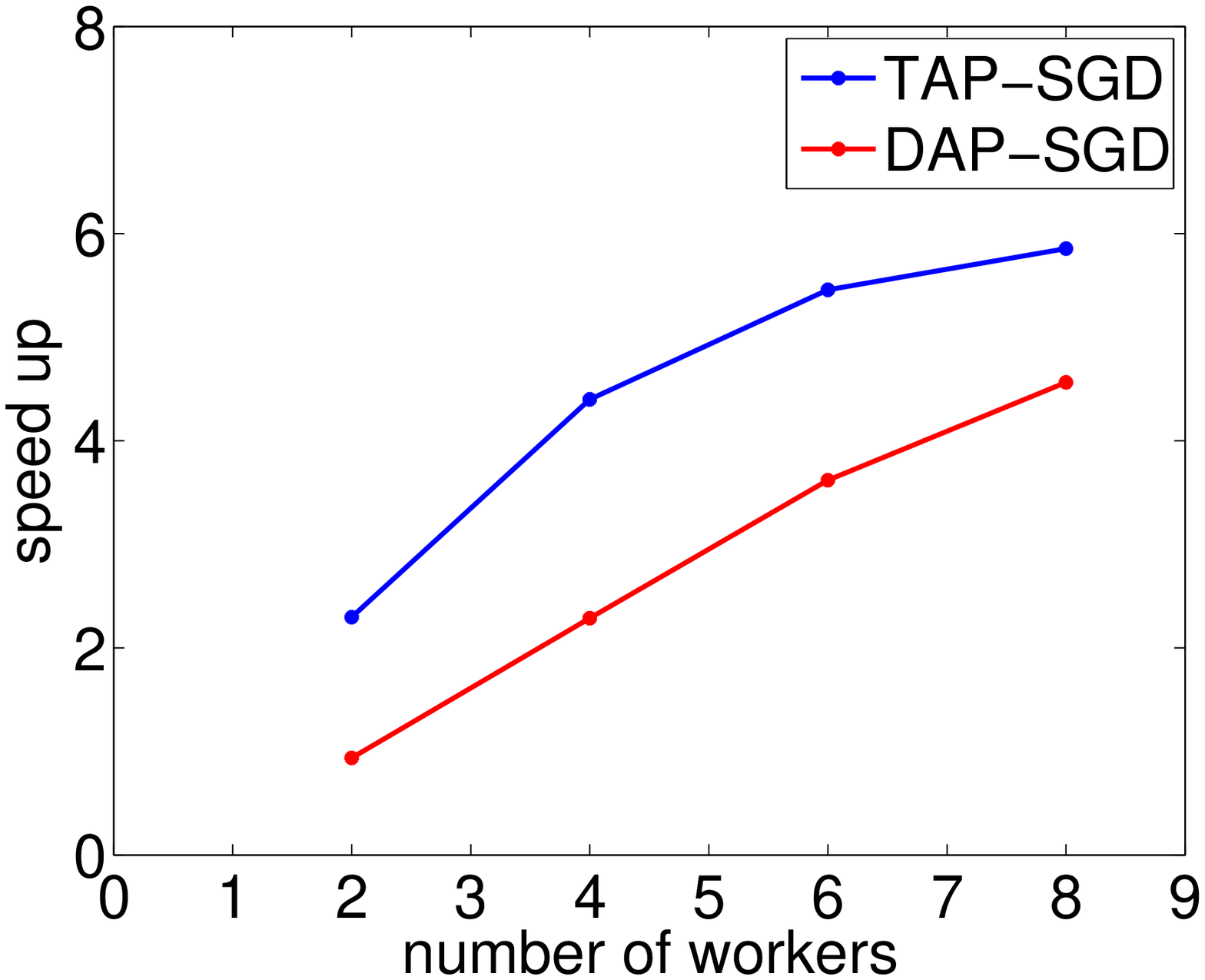}
} \subfigure[group lasso]{ \label{fig:exp1glasso}
\includegraphics[width=0.23\columnwidth]{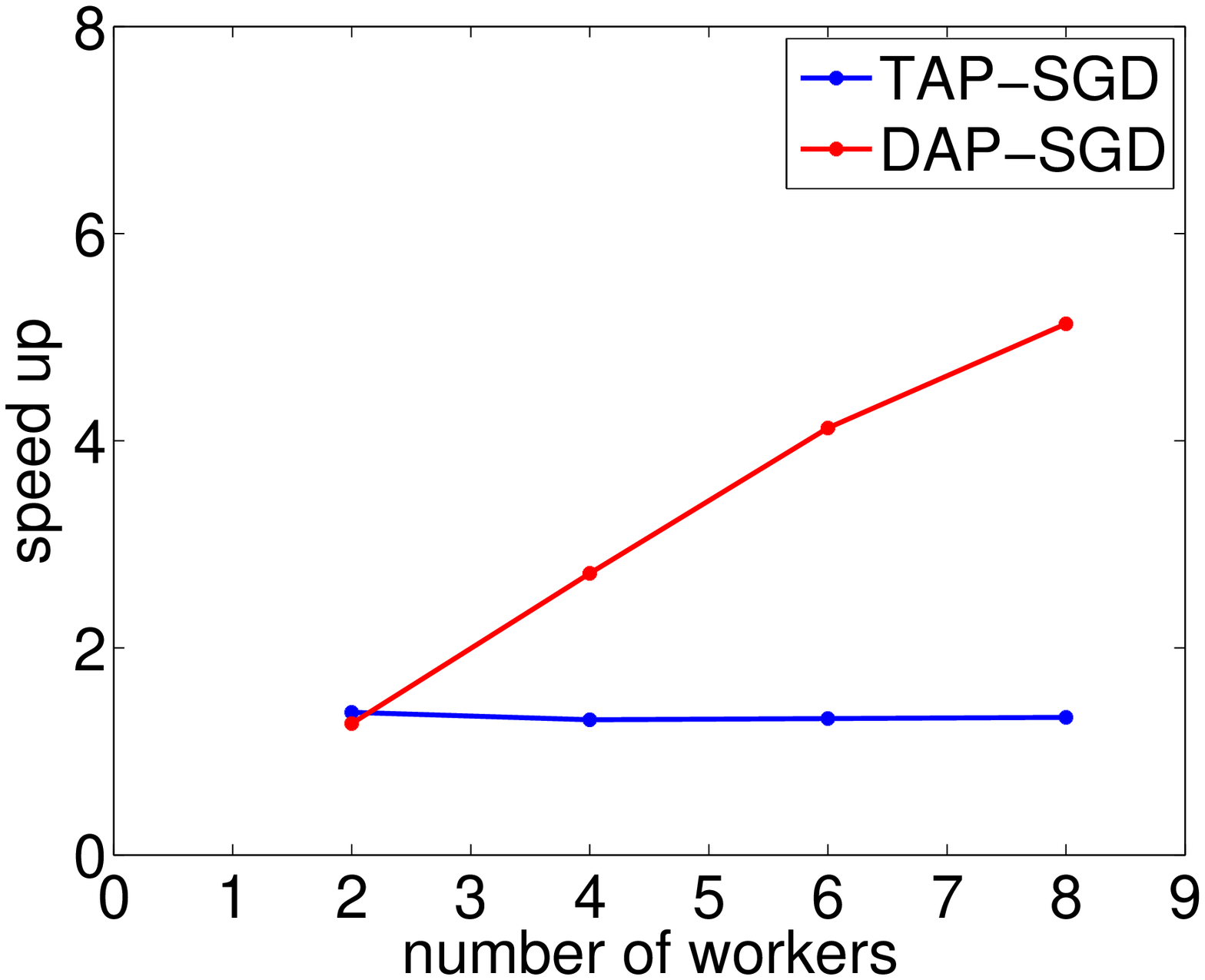}
} \subfigure[fused lasso]{ \label{fig:exp1sfusedlasso}
\includegraphics[width=0.23\columnwidth]{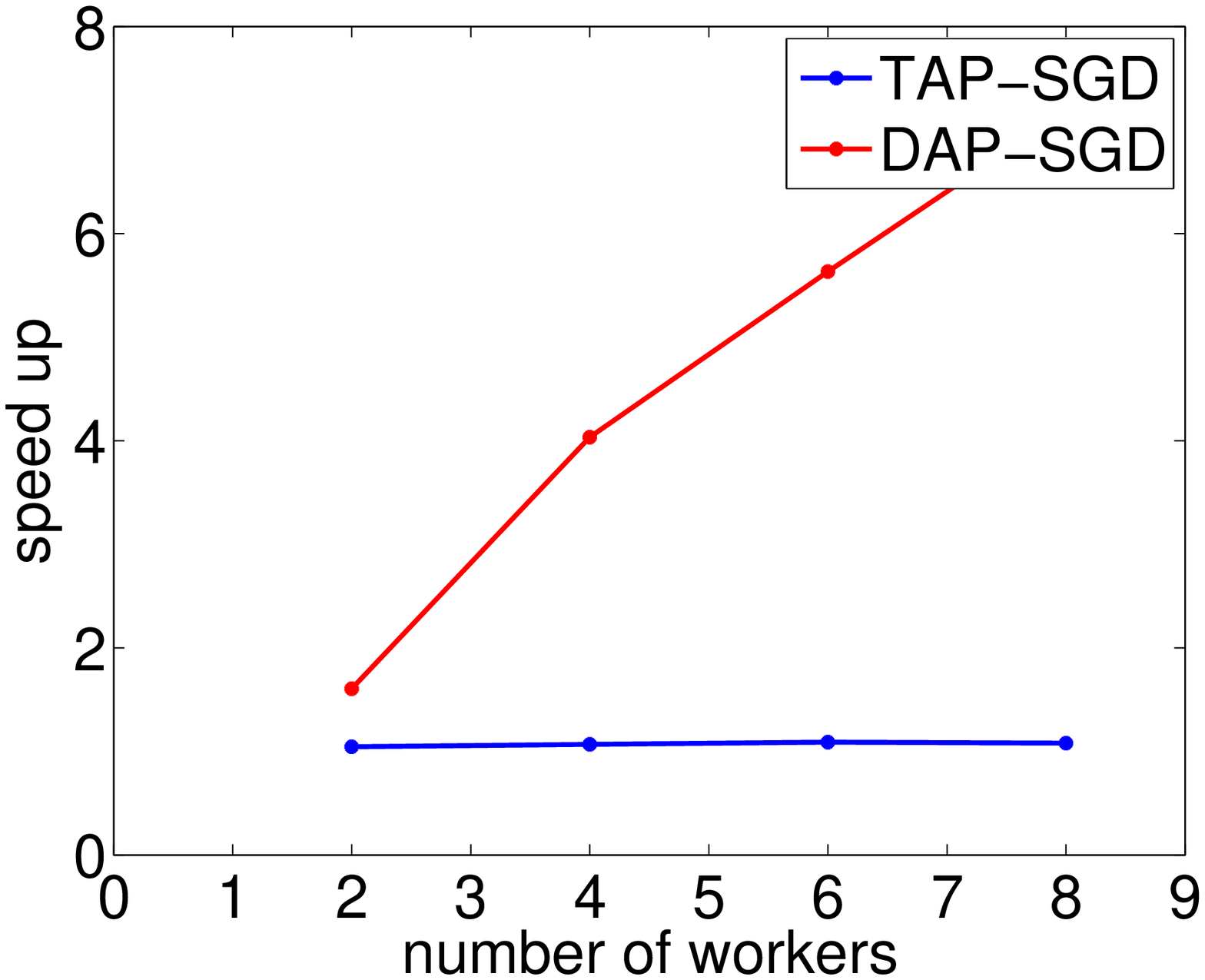}
} \subfigure[nuclear norm]{ \label{fig:exp1svt}
\includegraphics[width=0.23\columnwidth]{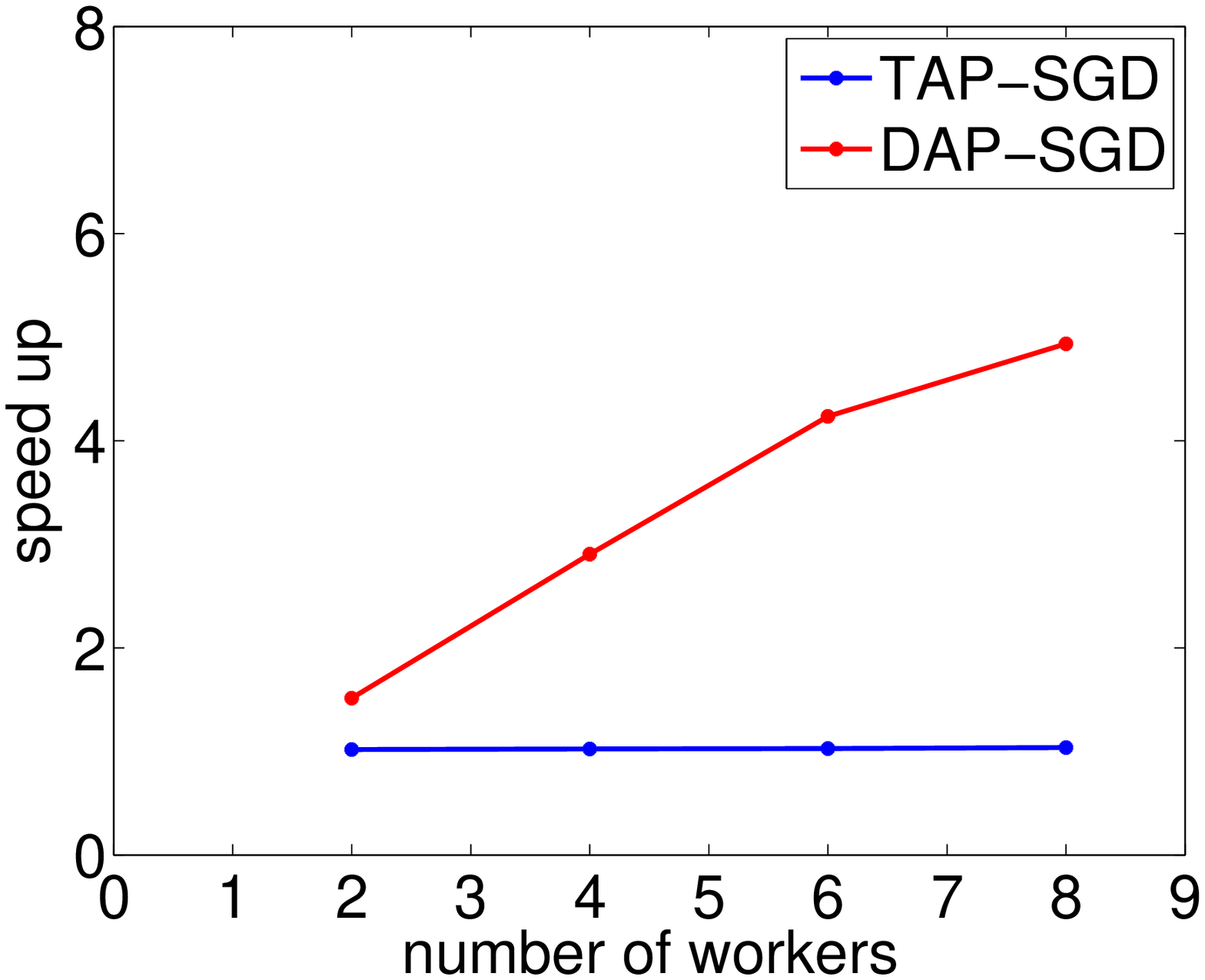}
} \caption{Speedup of TAP-SGD and DAP-SGD with 4 different
non-smooth regularization terms.}\label{fig:exp2}
\end{figure}

Results are summarized in Figures~\ref{fig:exp1} and
\ref{fig:exp2}. Figure \ref{fig:exp1} shows the comparison between
TAP-SGD and DAP-SGD regarding the running time and number of
iterations. As shown in the top row of Figure \ref{fig:exp1}, the
proposed DAP-GSD algorithm is slightly slower than TAP-SGD with
the $L_1$ regularized objective. The reason is that the proximal
operator of $L_1$ norm is element-wise and can be parallelized.
The decoupled update of DAP-SGD~(\ref{eqn:dapsgd}) involves more
operations in workers than the update of
TAP-SGD~(\ref{eqn:apsgd}), whose workers only need to evaluate the
gradients. Nevertheless, DAP-SGD is much faster than TAP-SGD with
group lasso, simplified fused lasso and nuclear norm regularized
objectives because the proximal operators of these norms are not
element-wise and hard to parallelize. As a consequence, evaluation
of the proximal operator in the master of TAP-SGD becomes the
computational bottleneck of the whole system and the performance
degrades significantly. In contrast, DAP-SGD allows each worker to
evaluate the proximal operator, which justifies our core idea of
decoupling the computation. Meanwhile, according to the bottom row
of Figure \ref{fig:exp1}, TAP-SGD and DAP-SGD perform similarly
regarding the number of iterations. The experimental results shown
in Figure~\ref{fig:exp1} validate that the decoupled operation in
DAP-SGD makes the algorithm more flexible and easier to
parallelize without affecting the precision of the algorithm.

Figure~\ref{fig:exp2} compares TAP-SGD and DAP-SGD in terms of the
speedup with different regularization terms. Obviously, DAP-SGD
can achieve significant speedup with the number of workers
increasing except for the $L_1$ regularized objective due to the
same reason discussed above. With group lasso, simplified fused
lasso and nuclear norm regularized objectives, TAP-SGD essentially
fails to speedup when the number of workers increases, which
indicates the computational bottleneck at the master for
evaluating the coupled proximal operator. Meanwhile, the
decoupling operation of DAP-SGD is effective to off-load the
computation to the workers and improves the parallelism in
asynchronous proximal stochastic gradient descent.

\section{Conclusion}

This paper proposes a novel decoupled asynchronous proximal
stochastic gradient descent~(DAP-SGD) algorithm for optimizing a
composite objective function. By off-loading computation from the
master to workers, the proposed DAP-SGD algorithm becomes easy to
parallelize. DAP-SGD is suitable for many master-worker
architectures, including single machine multi-core systems and
multi-machine systems. We further provide theoretical convergence
analyses for DAP-SGD, with both diminishing and fixed step-sizes.

\clearpage
\bibliographystyle{elsarticle-num}
\bibliography{ref}

\clearpage

\noindent\begin{Large}\textbf{\emph{Appendix} for
Make Workers Work Harder: Decoupled Asynchronous Proximal Stochastic Gradient Descent}\end{Large}

\begin{customthm}1
Suppose that the cost function of (1) satisfies the following
conditions: $f(\mathbf{x})$ is strongly convex with constant $\mu$
and $h(\mathbf{x})$ is convex; $f(\mathbf{x})$ is differentiable
and $\triangledown\! f(\mathbf{x})$ is Lipschitz continuous with
constant $L$; $\mathbb{E} \|\triangledown\!f_i(\mathbf{x}) -
\triangledown\!f(\mathbf{x})\|_2^2 \leq C_f$; $\|\partial
h(\mathbf{x})\|_2^2 \leq C_h$. Define the optimal solution of (1)
as $\mathbf{x}^*$. At time $t$, set the step-size of the DAP-SGD
recursion (8) as $\eta_t = O(1/t)$. Then the iterate generated by
(8) at time $T$, denoted by $\mathbf{x}_{T}$, satisfies
\begin{equation}
\begin{split}
\mathbb{E}\|\mathbf{x}_T - \mathbf{x}^*\|_2^2 \leq
O\left(\frac{\log T}{T}\right).
\end{split}
\end{equation}
\end{customthm}

\noindent\textbf{Proof of Theorem 1:} From the DAP-SGD update
$\mathbf{x}_{t+1} = \mathbf{x}_t + \mathbf{x}_{d(t)}' -
\mathbf{x}_{d(t)}$, we have
\begin{equation}\label{eqn:main}
\begin{split}
 &\mathbb{E} \|\mathbf{x}_{t+1} - \mathbf{x}^*\|_2^2 \\
=& \mathbb{E} \|\mathbf{x}_t - \mathbf{x}^* + \mathbf{x}_{d(t)}' - \mathbf{x}_{d(t)}\|^2 \\
=& \mathbb{E} \|\mathbf{x}_t - \mathbf{x}^*\|_2^2 + \mathbb{E} \|\mathbf{x}_{d(t)}' - \mathbf{x}_{d(t)}\|^2 + 2\mathbb{E}\left< \mathbf{x}_{d(t)}' - \mathbf{x}_{d(t)},\mathbf{x}_t - \mathbf{x}^* \right> \\
=& \mathbb{E} \|\mathbf{x}_t - \mathbf{x}^*\|_2^2 +
\underbrace{\mathbb{E} \|\mathbf{x}_{d(t)}' -
\mathbf{x}_{d(t)}\|^2 + 2\mathbb{E}\left< \mathbf{x}_{d(t)}' -
\mathbf{x}_{d(t)},\mathbf{x}_{d(t)} - \mathbf{x}^*
\right>}_{\mathrm{Q}_1} + 2\mathbb{E}\left< \mathbf{x}_{d(t)}' -
\mathbf{x}_{d(t)},\mathbf{x}_t - \mathbf{x}_{d(t)} \right>.
\end{split}
\end{equation}

Below we bound the value of $\mathrm{Q}_1$ from above.

Recalling the update of $\mathbf{x}'_{d(t)}$ in (8) of the paper, which is
\begin{equation}\label{eqn:ue1}
\begin{split}\mathbf{x}'_{d(t)} &= \text{Prox}_{\eta,h} (\mathbf{x}_{d(t)} - \eta_{d(t)} \triangledown\!f_{i_{d(t)}}(\mathbf{x}_{d(t)})) \\
&= \arg\!\min_{\mathbf{y}} \frac{1}{2\eta_{d(t)}}\|\mathbf{y} -
(\mathbf{x}_{d(t)} - \eta_{d(t)}
\triangledown\!f_{i_{d(t)}}(\mathbf{x}_{d(t)}))\|_2^2 +
h(\mathbf{y}),
\end{split}
\end{equation}
we have
\begin{equation}\label{eqn:prox_une}
\begin{split}
\frac{1}{\eta_{d(t)}}(\mathbf{x}_{d(t)} - \mathbf{x}'_{d(t)}) -
\triangledown f_{i_{d(t)}}(\mathbf{x}_{d(t)}) \in \partial
h(\mathbf{x}'_{d(t)}).
\end{split}
\end{equation}

Because $f(\mathbf{x})$ is convex (right now we do not need to use
its strong convexity) and $h(\mathbf{x})$ is also convex, we have
the following lower bound for the optimal value
\begin{equation}\label{ineqn:sc-0}
\begin{split}
P(\mathbf{x}^*) \triangleq & f(\mathbf{x}^*) + h(\mathbf{x}^*)  \\
\geq & f(\mathbf{x}_{d(t)}) + \left< \triangledown\!f(\mathbf{x}_{d(t)}),\mathbf{x}^* - \mathbf{x}_{d(t)} \right> + h(\mathbf{x}'_{d(t)}) + \left< \partial h(\mathbf{x}'_{d(t)}),\mathbf{x}^* - \mathbf{x}'_{d(t)} \right>. \\
\end{split}
\end{equation}
With a slight abuse of notation, here and thereafter $\partial
h(\mathbf{x}'_{d(t)})$ stands for any subgradient. Hence we
substitute the one given in (\ref{eqn:prox_une}) into
\eqref{ineqn:sc-0} and obtain
\begin{equation}\label{ineqn:sc}
\begin{split}
P(\mathbf{x}^*) \geq & f(\mathbf{x}_{d(t)}) + \left< \triangledown\!f(\mathbf{x}_{d(t)}),\mathbf{x}^* - \mathbf{x}_{d(t)} \right> \\
& + h(\mathbf{x}'_{d(t)}) + \left<
\frac{1}{\eta_{d(t)}}(\mathbf{x}_{d(t)} - \mathbf{x}'_{d(t)}) -
\triangledown\!f_{i_{d(t)}}(\mathbf{x}_{d(t)}),\mathbf{x}^* -
\mathbf{x}'_{d(t)} \right>.
\end{split}
\end{equation}

On the other hand, $\triangledown\!f(\mathbf{x})$ being Lipschitz
continuous with constant $L$ implies
\begin{equation}\label{ineqn:lcg}
\begin{split}
f(\mathbf{x}'_{d(t)}) \leq f(\mathbf{x}_{d(t)}) + \left<
\triangledown\!f(\mathbf{x}_{d(t)}),\mathbf{x}'_{d(t)} -
\mathbf{x}_{d(t)} \right> + \frac{L}{2}\|\mathbf{x}'_{d(t)} -
\mathbf{x}_{d(t)}\|_2^2.
\end{split}
\end{equation}

Substituting (\ref{ineqn:lcg}) into (\ref{ineqn:sc})
\begin{equation}\label{ineqn:sc-1}
\begin{split}
P(\mathbf{x}^*) \geq & f(\mathbf{x}'_{d(t)}) - \left<
\triangledown\!f(\mathbf{x}_{d(t)}),\mathbf{x}'_{d(t)} -
\mathbf{x}_{d(t)} \right> - \frac{L}{2}\|\mathbf{x}'_{d(t)} -
\mathbf{x}_{d(t)}\|_2^2 + \left< \triangledown\!f(\mathbf{x}_{d(t)}),\mathbf{x}^* - \mathbf{x}_{d(t)} \right> \\
&  + h(\mathbf{x}'_{d(t)}) + \left<
\frac{1}{\eta_{d(t)}}(\mathbf{x}_{d(t)} - \mathbf{x}'_{d(t)}) -
\triangledown\!f_{i_{d(t)}}(\mathbf{x}_{d(t)}),\mathbf{x}^* -
\mathbf{x}'_{d(t)} \right>.
\end{split}
\end{equation}
Noticing that by definition $P(\mathbf{x}'_{d(t)}) \triangleq
f(\mathbf{x}'_{d(t)}) + h(\mathbf{x}'_{d(t)})$ and reorganizing
the terms of \eqref{ineqn:sc-1}, we obtain
\begin{equation}\label{ineqn:t1-0}
\begin{split}
-[P(\mathbf{x}'_{d(t)}) - P(\mathbf{x}^*)] \geq & \left< \triangledown\!f(\mathbf{x}_{d(t)}) - \triangledown\!f_{i_{d(t)}}(\mathbf{x}_{d(t)}),\mathbf{x}^* - \mathbf{x}'_{d(t)} \right> + \frac{1}{\eta_{d(t)}}\left<\mathbf{x}_{d(t)} - \mathbf{x}'_{d(t)},\mathbf{x}^* - \mathbf{x}_{d(t)} \right> \\
& + \frac{1}{\eta_{d(t)}} \|\mathbf{x}_{d(t)} -
\mathbf{x}'_{d(t)}\|^2 - \frac{L}{2}\|\mathbf{x}_{d(t)} -
\mathbf{x}'_{d(t)}\|^2.
\end{split}
\end{equation}
Assuming that $\eta_t\leq 1/L$ for any $t$~(this assumption holds
according to the step-size rule given later), \eqref{ineqn:t1-0}
yields
\begin{equation}\label{ineqn:t1}
\begin{split}
-[P(\mathbf{x}'_{d(t)}) - P(\mathbf{x}^*)] \geq & \left< \triangledown\!f(\mathbf{x}_{d(t)}) - \triangledown\!f_{i_{d(t)}}(\mathbf{x}_{d(t)}),\mathbf{x}^* - \mathbf{x}'_{d(t)} \right> + \frac{1}{\eta_{d(t)}}\left< \mathbf{x}_{d(t)} - \mathbf{x}'_{d(t)},\mathbf{x}^* - \mathbf{x}_{d(t)} \right> \\
& + \frac{1}{2\eta_{d(t)}} \|\mathbf{x}_{d(t)} -
\mathbf{x}'_{d(t)}\|^2.
\end{split}
\end{equation}

Taking expectation on both sides of (\ref{ineqn:t1}) and
reorganizing terms, we have

\begin{equation}\label{ineqn:subt1}
\begin{split}
& -\mathbb{E}[P(\mathbf{x}'_{d(t)}) - P(\mathbf{x}^*)] + \underbrace{\mathbb{E}\left< \triangledown\!f_{i_{d(t)}}(\mathbf{x}_{d(t)}) - \triangledown\!f(\mathbf{x}_{d(t)}),\mathbf{x}^* - \mathbf{x}'_{d(t)} \right>}_{\mathrm{Q}_2} \\
\geq & \frac{1}{\eta_{d(t)}} \mathbb{E}\left< \mathbf{x}_{d(t)} -
\mathbf{x}'_{d(t)},\mathbf{x}^* - \mathbf{x}_{d(t)} \right> +
\frac{1}{2\eta_{d(t)}} \mathbb{E}\|\mathbf{x}_{d(t)} -
\mathbf{x}'_{d(t)}\|^2.
\end{split}
\end{equation}

Define $\hat{\mathbf{x}}'_{d(t)} \triangleq
\text{Prox}_{\eta,h}(\mathbf{x}_{d(t)} - \eta_{d(t)}
\triangledown\!f(\mathbf{x}_{d(t)}))$ as an approximation of
$\mathbf{x}'_{d(t)} \triangleq \text{Prox}_{\eta,h}
(\mathbf{x}_{d(t)} - \eta_{d(t)}
\triangledown\!f_{i_{d(t)}}(\mathbf{x}_{d(t)}))$. Because the
random variable $i_{d(t)}$ is independent with $\mathbf{x}^*$ and
$\hat{\mathbf{x}}'_{d(t)}$, while $\mathbb{E}
\left[\triangledown\!f_{i_{d(t)}}(\mathbf{x}_{d(t)})\right] =
\triangledown\!f(\mathbf{x}_{d(t)})$, it holds $\mathbb{E} \langle
\triangledown\!f_{i_{d(t)}}(\mathbf{x}_{d(t)}) -
\triangledown\!f(\mathbf{x}_{d(t)}),\mathbf{x}^* -
\hat{\mathbf{x}}'_{d(t)} \rangle = 0$. Hence, $\mathrm{Q}_2$ can
be upper bounded by
\begin{equation}\label{ineqn:non-exp-0}
\begin{split}
Q_2 = & \mathbb{E}\left< \triangledown\!f_{i_{d(t)}}(\mathbf{x}_{d(t)}) - \triangledown\!f(\mathbf{x}_{d(t)}),\mathbf{x}^* - \mathbf{x}'_{d(t)} \right> \\
=& \mathbb{E}\left< \triangledown\!f_{i_{d(t)}}(\mathbf{x}_{d(t)}) - \triangledown\!f(\mathbf{x}_{d(t)}),\mathbf{x}^* - \hat{\mathbf{x}}'_{d(t)} \right>+ \mathbb{E}\left< \triangledown\!f_{i_{d(t)}}(\mathbf{x}_{d(t)}) - \triangledown\!f(\mathbf{x}_{d(t)}), \hat{\mathbf{x}}'_{d(t)} - \mathbf{x}'_{d(t)} \right> \\
=& \mathbb{E}\left< \triangledown\!f_{i_{d(t)}}(\mathbf{x}_{d(t)}) - \triangledown\!f(\mathbf{x}_{d(t)}), \hat{\mathbf{x}}'_{d(t)} - \mathbf{x}'_{d(t)} \right> \\
\leq & \mathbb{E}
\left(\|\triangledown\!f_{i_{d(t)}}(\mathbf{x}_{d(t)}) -
\triangledown\!f(\mathbf{x}_{d(t)})\|_2\|\hat{\mathbf{x}}'_{d(t)}
- \mathbf{x}'_{d(t)}\|_2 \right),
\end{split}
\end{equation}
where the last inequality comes from the Cauchy-Schwarz inequality.
Further, the non-expansive property of proximal operators [8]
implies
\begin{equation}\label{ineqn:non-exp-x}
\begin{split}
\|\hat{\mathbf{x}}'_{d(t)} - \mathbf{x}'_{d(t)}\|_2 = & \|
\text{Prox}_{\eta,h}(\mathbf{x}_{d(t)} - \eta_{d(t)}
\triangledown\!f(\mathbf{x}_{d(t)})) - \text{Prox}_{\eta,h}
(\mathbf{x}_{d(t)} - \eta_{d(t)}
\triangledown\!f_{i_{d(t)}}(\mathbf{x}_{d(t)})) \|_2 \\
\leq& \eta_{d(t)} \|\triangledown\!f_{i_{d(t)}}(\mathbf{x}_{d(t)})
- \triangledown\!f(\mathbf{x}_{d(t)})\|_2 .
\end{split}
\end{equation}
Combining \eqref{ineqn:non-exp-0} and \eqref{ineqn:non-exp-x}
yields an upper bound of $Q_2$ as
\begin{equation}\label{ineqn:non-exp}
\begin{split}
Q_2  \leq  \eta_{d(t)}\mathbb{E}
\|\triangledown\!f_{i_{d(t)}}(\mathbf{x}_{d(t)}) -
\triangledown\!f(\mathbf{x}_{d(t)})\|_2^2\leq  \eta_{d(t)} C_f,
\end{split}
\end{equation}
where the last inequality is due to the assumption of bounded variance
$\mathbb{E} \|\triangledown\!f_i(\mathbf{x}) -
\triangledown\!f(\mathbf{x})\|_2^2 \leq C_f$.

Substituting (\ref{ineqn:non-exp}) into (\ref{ineqn:subt1}), we
have
\begin{equation}\label{ineqn:lastt1}
\begin{split}
& -\mathbb{E}[P(\mathbf{x}'_{d(t)}) - P(\mathbf{x}^*)] + \eta_{d(t)} C_f \\
\geq & \frac{1}{\eta_{d(t)}} \mathbb{E}\left< \mathbf{x}_{d(t)} -
\mathbf{x}'_{d(t)},\mathbf{x}^* - \mathbf{x}_{d(t)} \right> +
\frac{1}{2\eta_{d(t)}} \mathbb{E}\|\mathbf{x}_{d(t)} -
\mathbf{x}'_{d(t)}\|^2.
\end{split}
\end{equation}

Now we end up with an upper bound of $\mathrm{Q}_1$ as
\begin{equation}\label{eqn:svrg-0}
\begin{split}
\mathrm{Q}_1 \triangleq & \mathbb{E} \|\mathbf{x}_{d(t)}' - \mathbf{x}_{d(t)}\|^2 + 2\mathbb{E}\left< \mathbf{x}_{d(t)}' - \mathbf{x}_{d(t)},\mathbf{x}_{d(t)} - \mathbf{x}^* \right> \\
\leq & - 2\eta_{d(t)}\mathbb{E}[P(\mathbf{x}'_{d(t)}) -
P(\mathbf{x}^*)] + 2\eta_{d(t)}^2 C_f.
\end{split}
\end{equation}
Therefore
\begin{equation}\label{eqn:svrg}
\begin{split}
\mathrm{Q}_1 \leq &  - 2\eta_{d(t)}\mathbb{E}[P(\mathbf{x}_t) -
P(\mathbf{x}^*)] - 2\eta_{d(t)}\mathbb{E}[P(\mathbf{x}'_{d(t)}) -
P(\mathbf{x}_t)] + 2\eta_{d(t)}^2 C_f. \\
\leq &  - \mu\eta_{d(t)}\mathbb{E}\|\mathbf{x}_t -
\mathbf{x}^*\|_2^2 - 2\eta_{d(t)}\mathbb{E}[P(\mathbf{x}'_{d(t)})
- P(\mathbf{x}_t)] + 2\eta_{d(t)}^2 C_f.
\end{split}
\end{equation}
The second line comes from the inequality
\begin{equation}\label{eqn:svrg-1}
\begin{split}
P(\mathbf{x}_t) - P(\mathbf{x}^*) \geq \frac{\mu}{2}\|\mathbf{x}_t
- \mathbf{x}^*\|_2^2,
\end{split}
\end{equation}
which is due to the facts that $\mathbf{x}^*$ is the optimal
solution of $P(\mathbf{x})=f(\mathbf{x})+h(\mathbf{x})$,
$f(\mathbf{x})$ is strongly convex with constant $\mu$, and
$h(\mathbf{x})$ is convex.

Substituting (\ref{eqn:svrg}) into (\ref{eqn:main}), we have
\begin{equation}\label{eqn:main2}
\begin{split}
\mathbb{E} \|\mathbf{x}_{t+1} - \mathbf{x}^*\|_2^2 \leq & (1-\mu \eta_{d(t)}) \mathbb{E} \|\mathbf{x}_t - \mathbf{x}^*\|_2^2 + 2 \eta_{d(t)} \underbrace{\mathbb{E} [P(\mathbf{x}_{d(t)}) - P(\mathbf{x}'_{d(t)})]}_{\mathrm{Q}_3} \\
& + 2 \eta_{d(t)} \sum_{p = 1}^{t - d(t)} \underbrace{\mathbb{E}
[P(\mathbf{x}_{t-p+1}) - P(\mathbf{x}_{t-p})]}_{\mathrm{Q}_4} +
2\eta_{d(t)}^2 C_f + 2\underbrace{\mathbb{E}\left<
\mathbf{x}_{d(t)}' - \mathbf{x}_{d(t)},\mathbf{x}_t -
\mathbf{x}_{d(t)} \right>}_{\mathrm{Q}_5}.
\end{split}
\end{equation}
We proceed to bound the terms $\mathrm{Q}_3$, $\mathrm{Q}_4$, and
$\mathrm{Q}_5$.

Because $f(\mathbf{x})$ and $h(\mathbf{x})$ are convex as well as
the norm of $\partial h(\mathbf{x})$ is bounded, we have the
following basic inequality
\begin{equation}\label{eqn:subb}
\begin{split}
P(\mathbf{x}) - P(\mathbf{y}) = & f(\mathbf{x}) - f(\mathbf{y}) + h(\mathbf{x}) - h(\mathbf{y}) \\
\leq & \left< \triangledown\!f(\mathbf{x}), \mathbf{x} - \mathbf{y} \right> +  \left< \partial h(\mathbf{x}),\mathbf{x} - \mathbf{y} \right> \\
\leq & \|\triangledown\!f(\mathbf{x})\|_2\|\mathbf{x} - \mathbf{y}\|_2 +  \|\partial h(\mathbf{x})\|_2\|\mathbf{x} - \mathbf{y}\|_2 \\
\leq & \|\triangledown\!f(\mathbf{x})\|_2 \|\mathbf{x} - \mathbf{y}\|_2 + \sqrt{C_h} \|\mathbf{x} - \mathbf{y}\|_2 \\
= & (\|\triangledown\!f(\mathbf{x})\|_2 + \sqrt{C_h}) \|\mathbf{x}
- \mathbf{y}\|_2.
\end{split}
\end{equation}
In \eqref{eqn:subb}, the second line comes from the convexity of
$f(\mathbf{x})$ and $h(\mathbf{x})$, while the third line comes
from the Cauchy-Schwarz inequality. Replacing $\mathbf{x}$ by
$\mathbf{x}_{d(t)}$ and $\mathbf{y}$ by $\mathbf{x}'_{d(t)}$ in
(\ref{eqn:subb}), we have
\begin{equation}\label{eqn:t3-0}
\begin{split}
\mathrm{Q}_3 = \mathbb{E} \left[ P(\mathbf{x}_{d(t)}) -
P(\mathbf{x}'_{d(t)}) \right] \leq
\mathbb{E}\left[(\|\triangledown f(\mathbf{x}_{d(t)})\|_2 +
\sqrt{C_h})\|\mathbf{x}_{d(t)} - \mathbf{x}'_{d(t)}\|_2\right].
\end{split}
\end{equation}
Applying the expression of $\mathbf{x}_{d(t)} -
\mathbf{x}'_{d(t)}$ in \eqref{eqn:prox_une} into \eqref{eqn:t3-0}
yields
\begin{equation}\label{eqn:t3-1}
\begin{split}
\mathrm{Q}_3 \leq & \eta_{d(t)} \mathbb{E}\left[ (\|\triangledown\!f(\mathbf{x}_{d(t)})\|_2 + \sqrt{C_h})\|\triangledown\!f_{i_{d(t)}}(\mathbf{x}_{d(t)}) + \partial h(\mathbf{x}'_{d(t)})\|_2\right] \\
\leq & \frac{1}{2} \eta_{d(t)}
\mathbb{E}\|\triangledown\!f(\mathbf{x}_{d(t)})\|_2^2 +
\frac{1}{2}  \eta_{d(t)} C_h + \eta_{d(t)}\mathbb{E}
\|\triangledown\!f_{i_{d(t)}}(\mathbf{x}_{d(t)}) + \partial
h(\mathbf{x}'_{d(t)})\|_2^2.
\end{split}
\end{equation}
Due to the inequalities
\begin{equation}\label{eqn:t3-xx}
\begin{split}
\frac{1}{2}\|\triangledown\!f(\mathbf{x}_{d(t)})\|_2^2 \leq
\|\triangledown\!f(\mathbf{x}_{d(t)}) -
\triangledown\!f(\mathbf{x}^*)\|_2^2 + \|
\triangledown\!f(\mathbf{x}^*)\|_2^2,
\end{split}
\end{equation}
and
\begin{equation}\label{eqn:t3-yy}
\begin{split}
& \|\triangledown\!f_{i_{d(t)}}(\mathbf{x}_{d(t)}) + \partial
h(\mathbf{x}'_{d(t)})\|_2^2 \\
\leq &2 \|\triangledown\!f_{i_{d(t)}}(\mathbf{x}_{d(t)})\|_2^2 +
2\|\partial h(\mathbf{x}'_{d(t)})\|_2^2 \\
\leq & 4\|\triangledown\!f_{i_{d(t)}}(\mathbf{x}_{d(t)}) -
\triangledown\!f(\mathbf{x}_{d(t)})\|_2^2 +
4\|\triangledown\!f(\mathbf{x}_{d(t)})\|_2^2 +
2\|\partial h(\mathbf{x}'_{d(t)})\|_2^2 \\
\leq & 4\|\triangledown\!f_{i_{d(t)}}(\mathbf{x}_{d(t)}) -
\triangledown\!f(\mathbf{x}_{d(t)})\|_2^2 +
8\|\triangledown\!f(\mathbf{x}_{d(t)}) - \triangledown
f(\mathbf{x}^*)\|_2^2 + 8\|\triangledown\!f(\mathbf{x}^*)\|_2^2 +
2\|\partial h(\mathbf{x}'_{d(t)})\|_2^2,
\end{split}
\end{equation}
\eqref{eqn:t3-1} turns to
\begin{equation}\label{eqn:t3-2}
\begin{split}
\mathrm{Q}_3 \leq & 9\eta_{d(t)} \mathbb{E}\|\triangledown\!f(\mathbf{x}_{d(t)}) - \triangledown\!f(\mathbf{x}^*)\|_2^2 + 9\eta_{d(t)} \mathbb{E}\| \triangledown\!f(\mathbf{x}^*)\|_2^2 + 4\eta_{d(t)} \mathbb{E}\|\triangledown\!f_{i_{d(t)}}(\mathbf{x}_{d(t)}) - \triangledown\!f(\mathbf{x}_{d(t)})\|_2^2 \\
& + 2\eta_{d(t)}\mathbb{E}\|\partial h(\mathbf{x}'_{d(t)})\|_2^2 + \frac{1}{2} \eta_{d(t)}C_h. \\
\end{split}
\end{equation}
Considering Lipschitz continuity of $\triangledown\!f(\mathbf{x})$,
$\|\triangledown\!f(\mathbf{x}^*)\|_2^2 \leq C_h$ from Corollary
1, $\mathbb{E} \|\triangledown\!f_i(\mathbf{x}) -
\triangledown\!f(\mathbf{x})\|_2^2 \leq C_f$, as well as
$\|\partial h(\mathbf{x})\|_2^2 \leq C_h$, \eqref{eqn:t3-2}
further turns to
\begin{equation}\label{eqn:t3}
\begin{split}
\mathrm{Q}_3 \leq 9 \eta_{d(t)} L^2 \mathbb{E}\|\mathbf{x}_{d(t)}
- \mathbf{x}^*\|_2^2 + 4\eta_{d(t)} C_f + \frac{23}{2} \eta_{d(t)}
C_h.
\end{split}
\end{equation}

Similar to the derivation of (\ref{eqn:t3}), we have
\begin{equation}\label{eqn:t4-0}
\begin{split}
\mathrm{Q}_4 =& \mathbb{E}\left[ P(\mathbf{x}_{t-p+1}) - P(\mathbf{x}_{t-p}) \right] \\
\leq & \mathbb{E}\left[ (\|\triangledown f(\mathbf{x}_{t-p+1})\|_2 + \sqrt{C_h}) \|\mathbf{x}_{t-p+1} - \mathbf{x}_{t-p}\|_2 \right] \\
\leq & \eta_{d(t-p)}\mathbb{E}
\left[(\|\triangledown\!f(\mathbf{x}_{t-p+1})\|_2 + \sqrt{C_h})
\|\triangledown\!f_{i_{d(t-p)}}(\mathbf{x}_{d(t-p)}) + \partial
h(\mathbf{x}'_{d(t-p)})\|_2 \right] \\
\leq & \frac{1}{2} \eta_{d(t-p)}
\mathbb{E}\|\triangledown\!f(\mathbf{x}_{t-p+1})\|_2^2 +
\frac{1}{2} \eta_{d(t-p)} C_h + \eta_{d(t-p)}
\mathbb{E}\|\triangledown\!f_{i_{d(t-p)}}(\mathbf{x}_{d(t-p)}) +
\partial h(\mathbf{x}'_{d(t-p)})\|_2^2.
\end{split}
\end{equation}
Using the inequalities (see \eqref{eqn:t3-xx} and
\eqref{eqn:t3-yy})
\begin{equation}\label{eqn:t4-1}
\begin{split}
\frac{1}{2} \|\triangledown\!f(\mathbf{x}_{t-p+1})\|_2^2 \leq
\|\triangledown\!f(\mathbf{x}_{t-p+1}) -
\triangledown\!f(\mathbf{x}^*)\|_2^2 +
\|\triangledown\!f(\mathbf{x}^*)\|_2^2,
\end{split}
\end{equation}
and
\begin{equation}\label{eqn:t4-2}
\begin{split}
& \|\triangledown\!f_{i_{d(t-p)}}(\mathbf{x}_{d(t-p)}) +
\partial h(\mathbf{x}'_{d(t-p)})\|_2^2 \\
\leq & 4\|\triangledown\!f_{i_{d(t-p)}}(\mathbf{x}_{d(t-p)}) -
\triangledown\!f(\mathbf{x}_{d(t-p)})\|_2^2 +
8\|\triangledown\!f(\mathbf{x}_{d(t-p)}) -
\triangledown\!f(\mathbf{x}^*)\|_2^2
+8\|\triangledown\!f(\mathbf{x}^*)\|_2^2 + 2\|\partial
h(\mathbf{x}'_{d(t-p)})\|_2^2,
\end{split}
\end{equation}
\eqref{eqn:t4-0} yields
\begin{equation}\label{eqn:t4}
\begin{split}
\mathrm{Q}_4 \leq &
\eta_{d(t-p)}\mathbb{E}\|\triangledown\!f(\mathbf{x}_{t-p+1}) -
\triangledown\!f(\mathbf{x}^*)\|_2^2 + 9
\eta_{d(t-p)}\mathbb{E}\|\triangledown\!f(\mathbf{x}^*)\|^2 + 8
\eta_{d(t-p)}\mathbb{E}\|\triangledown\!f(\mathbf{x}_{t-p}) -
\triangledown\!f(\mathbf{x}^*)\|_2^2 \\
& + 4\eta_{d(t-p)} \mathbb{E} \|\triangledown\!f_{i_{d(t-p)}}(\mathbf{x}_{d(t-p)}) - \triangledown\!f(\mathbf{x}_{d(t-p)})\|_2^2 + 2\eta_{d(t-p)}\mathbb{E}\|\partial h(\mathbf{x}'_{d(t-p)})\|_2^2 + \frac{1}{2} \eta_{d(t-p)} C_h \\
\leq & \eta_{d(t-p)}L^2 \mathbb{E}\|\mathbf{x}_{t-p+1} -
\mathbf{x}^*\|_2^2 + 8\eta_{d(t-p)} L^2
\mathbb{E}\|\mathbf{x}_{d(t-p)} - \mathbf{x}^*\|_2^2 +
4\eta_{d(t-p)} C_f + \frac{23}{2} \eta_{d(t-p)} C_h.
\end{split}
\end{equation}
Again, the last line of \eqref{eqn:t4} utilizes Lipschitz
continuity of $\triangledown\!f(x)$,
$\|\triangledown\!f(\mathbf{x}^*)\|_2^2 \leq C_h$ from Corollary
1, $\mathbb{E} \|\triangledown\!f_i(\mathbf{x}) -
\triangledown\!f(\mathbf{x})\|_2^2 \leq C_f$, as well as
$\|\partial h(\mathbf{x})\|_2^2 \leq C_h$.

For the term $Q_5$, we use the Cauchy-Schwarz inequality followed
by the substitution of \eqref{eqn:prox_une} and get
\begin{equation}\label{eqn:t5-0}
\begin{split}
\mathrm{Q}_5 =& \mathbb{E} \left< \mathbf{x}_{d(t)}' - \mathbf{x}_{d(t)}, \mathbf{x}_t - \mathbf{x}_{d(t)} \right> \\
\leq & \mathbb{E} \left( \|\mathbf{x}_{d(t)}' - \mathbf{x}_{d(t)}\|_2 \|\mathbf{x}_t - \mathbf{x}_{d(t)}\|_2 \right) \\
\leq & \eta_{d(t)} \mathbb{E} \left( \|\triangledown\!f_{i_{d(t)}}(\mathbf{x}_{d(t)}) + \partial h(\mathbf{x}'_{d(t)})\|_2 \|\mathbf{x}_t - \mathbf{x}_{d(t)}\|_2 \right). \\
\end{split}
\end{equation}
Further relaxing \eqref{eqn:t5-0} by the triangle inequality
yields
\begin{equation}\label{eqn:t5-1}
\begin{split}
\mathrm{Q}_5 \leq \eta_{d(t)} \sum_{p = 1}^{t - d(t)} \mathbb{E}
\left( \|\triangledown\!f_{i_{d(t)}}(\mathbf{x}_{d(t)}) +
\partial h(\mathbf{x}'_{d(t)})\|_2 \|\mathbf{x}_{t-p+1} -
\mathbf{x}_{t-p}\|_2 \right).
\end{split}
\end{equation}
Since the maximum delay is $\tau$, we have
\begin{equation}\label{eqn:t5-2}
\begin{split}
\mathrm{Q}_5 \leq \eta_{d(t)} \sum_{p = 1}^{\tau} \mathbb{E} \left( \|\triangledown\!f_{i_{d(t)}}(\mathbf{x}_{d(t)}) + \partial h(\mathbf{x}'_{d(t)})\|_2 \|\mathbf{x}_{t-p+1} - \mathbf{x}_{t-p}\|_2 \right). \\
\end{split}
\end{equation}
Noticing the relations $\mathbf{x}_{t-p+1} - \mathbf{x}_{t-p} =
\mathbf{x}'_{d(t-p)} - \mathbf{x}_{d(t-p)}$ from the DAP-SGD
recursion and $\mathbf{x}'_{d(t-p)} - \mathbf{x}_{d(t-p)} =
\eta_{d(t)} (\triangledown\!f(\mathbf{x}_{d(t-p)}) + \partial
h(\mathbf{x}'_{d(t-p)}))$ from \eqref{eqn:prox_une},
\eqref{eqn:t5-2} leads to
\begin{equation}\label{eqn:t5-2}
\begin{split}
\mathrm{Q}_5 \leq \eta_{d(t)} \sum_{p = 1}^{\tau} \eta_{d(t-p)}
\mathbb{E} \left( \|\triangledown\!f_{i_{d(t)}}(\mathbf{x}_{d(t)})
+ \partial h(\mathbf{x}'_{d(t)})\|_2
\|\triangledown\!f(\mathbf{x}_{d(t-p)}) + \partial
h(\mathbf{x}'_{d(t-p)})\|_2 \right).
\end{split}
\end{equation}
Following the similar routines as those in (\ref{eqn:t3})\ and
(\ref{eqn:t4}), eventually we reach
\begin{equation}\label{eqn:t5}
\begin{split}
\mathrm{Q}_5 \leq & 4 \eta_{d(t)}L^2 \sum_{p = 1}^{\tau} \eta_{d(t-p)}\mathbb{E} \|\mathbf{x}_{d(t)} - \mathbf{x}^*\|_2^2 + 4 \eta_{d(t)}L^2\sum_{p = 1}^{\tau} \eta_{d(t-p)} \mathbb{E}\|\mathbf{x}_{d(t-p)} - \mathbf{x}^*\|_2^2 \\
& + 4 \eta_{d(t)}\sum_{p = 1}^{\tau} \eta_{d(t-p)} C_f +
10\eta_{d(t)} \sum_{p = 1}^{\tau} \eta_{d(t-p)} C_h
\end{split}
\end{equation}

Substituting (\ref{eqn:t3}), (\ref{eqn:t4}) and (\ref{eqn:t5})
into (\ref{eqn:main2}), we have
\begin{equation}\label{eqn:main3}
\begin{split}
\mathbb{E} \|\mathbf{x}_{t+1} - \mathbf{x}^*\|_2^2 \leq & \left(1-\mu \eta_{d(t)} \right) \mathbb{E} \|\mathbf{x}_t - \mathbf{x}^*\|_2^2 + \left(8 \eta_{d(t)}L^2 \sum_{p = 1}^{\tau} \eta_{d(t-p)} + 18 \eta_{d(t)}^2 L^2 \right) \mathbb{E}\|\mathbf{x}_{d(t)} - \mathbf{x}^*\|_2^2 \\
& + 2 \eta_{d(t)} L^2\sum_{p = 1}^{\tau} \eta_{d(t-p)} \mathbb{E}\|\mathbf{x}_{t-p+1} - \mathbf{x}^*\|_2^2 + 24\eta_{d(t)} L^2\sum_{p = 1}^{\tau} \eta_{d(t-p)} \mathbb{E}\|\mathbf{x}_{d(t-p)} - \mathbf{x}^*\|_2^2 \\
& + \left( 16 \eta_{d(t)}\sum_{p = 1}^{\tau} \eta_{d(t-p)} +
8\eta_{d(t)}^2 \right) C_f + \left( 43\eta_{d(t)}\sum_{p =
1}^{\tau} \eta_{d(t-p)} + 23 \eta_{d(t)}^2 \right) C_h.
\end{split}
\end{equation}

Define the step-size rule
\begin{equation}\label{eqn:stepsizerule}
\begin{split}
\eta_t = \frac{1}{\mu (t + 1) + u} = O\left(\frac{1}{t}\right),
\end{split}
\end{equation}
where $u$ is a positive constant satisfying:
\begin{itemize}
\item $u > (2\tau - 1)\mu$ such that $\eta_t \leq \eta_{d(t)}$;
\item $u$ is large enough such that $\min(\mu/(4 C_1 \tau),
1/L)\geq \eta_t$, where $C_1$ is a constant we give below.
\end{itemize}
Define two constants
$$C_1 = \left(2L^2 \frac{\mu + u}{\mu + u - 2\mu\tau} + 48\tau L^2 +
8\tau L^2 \frac{\mu + u}{\mu + u - 2\mu\tau}\right) \frac{\mu +
u}{\mu + u - \mu\tau} + 18L^2, $$
and
$$C_2 = \left[(16 \tau + 8) C_f + (43 \tau + 23) C_h
\right]\frac{(\mu+u)^2}{{(\mu + u - 2\mu\tau)}^2}.$$
Though not straightforward, we can show that under the step-size
rule given by \eqref{eqn:stepsizerule}, \eqref{eqn:main3} yields
\begin{equation}\label{eqn:thm1basic-0}
\begin{split}
\mathbb{E} \|\mathbf{x}_{t+1} - \mathbf{x}^*\|_2^2 \leq (1-\mu
\eta_t) \mathbb{E} \|\mathbf{x}_t - \mathbf{x}^*\|_2^2 + C_1
\sum_{p = 0}^{2\tau} \eta_{t-p}^2 \mathbb{E} \|\mathbf{x}_{t-p} -
\mathbf{x}^*\|_2^2 + C_2 \eta_t^2.
\end{split}
\end{equation}

For the ease of presentation, we define $a_t = \mathbb{E}
\|\mathbf{x}_t - \mathbf{x}^*\|_2^2$ and will analyze its rate.
Rewrite \eqref{eqn:thm1basic-0} to
\begin{equation}\label{eqn:thm1basic}
\begin{split}
a_{t+1} \leq (1-\mu \eta_t) a_t + C_1 \sum_{p = 0}^{2\tau}
\eta_{t-p}^2 a_{t-p} + C_2 \eta_t^2.
\end{split}
\end{equation}

Applying telescopic cancellation to (\ref{eqn:thm1basic}) from $t
= 0$ to $t = T-1$ yields
\begin{equation} \label{eqn:aaaa-1}
\begin{split}
a_T \leq & a_0 - \sum_{t=0}^{T-1} \mu \eta_t a_t + C_1
\sum_{t=0}^{T-1} \sum_{p = 0}^{2\tau} \eta_{t-p}^2 a_{t-p} + C_2
\sum_{t=0}^{T-1} \eta_t^2 \\
\leq & a_0 - \sum_{t=0}^{T-1} (\mu\eta_t - 2C_1 \eta_t^2\tau)a_t +
C_2 O(1).
\end{split}
\end{equation}
As we can verify, $\mu/(4 C_1 \tau)\geq \eta_t$, meaning that
\begin{equation} \label{eqn:aaaa-2}
\begin{split}
\sum_{t=0}^{T-1} (\mu\eta_t - 2C_1 \eta_t^2\tau)a_t \geq
\frac{1}{2}\sum_{t=0}^{T-1} \mu\eta_ta_t.
\end{split}
\end{equation}
Combining \eqref{eqn:aaaa-1} and \eqref{eqn:aaaa-2}, we have
\begin{equation}\label{ineqn:aaaa-3}
\begin{split}
\frac{1}{2}\sum_{t=0}^{T-1} \mu\eta_ta_t \leq a_0 - a_{T} + C_2
O(1),
\end{split}
\end{equation}
which, along with the step-size rule (\ref{eqn:stepsizerule}), implies that
\begin{equation}\label{ineqn:atbnd}
\begin{split}
\sum_{t=0}^{T-1} \frac{1}{\mu(t+1) + u} a_t \leq \frac{2}{\mu}
(a_0 + C_2 O(1))
\end{split}
\end{equation}

Further define $C_3 = u /(u - \mu\tau)$ such that $$
\frac{\mu(t+1) + u}{(\mu(t-p+1) + u)^2}\leq \frac{C_3}{\mu(t-p+1)
+ u}.$$ Substituting the step-size rule (\ref{eqn:stepsizerule}) into
(\ref{eqn:thm1basic}), we have
\begin{equation}\label{eqn:thm1basic2-0}
\begin{split}
a_{t+1} \leq & \left(1-\frac{\mu}{\mu (t+1) + u}\right) a_t + C_1
\sum_{p=0}^{2\tau} \frac{1}{(\mu(t-p+1) + u)^2} a_{t-p} +
\frac{1}{(\mu (t+1) + u)^2} C_2,
\end{split}
\end{equation}
and consequently
\begin{equation}\label{eqn:thm1basic2}
\begin{split}
(\mu (t+1) + u) a_{t+1} \leq & (\mu t + u) a_t + C_1 \sum_{p=0}^{2\tau} \frac{\mu (t+1) + u}{(\mu(t-p+1) + u)^2} a_{t-p} + \frac{1}{\mu (t+1) + u} C_2 \\
\leq & (\mu t + u) a_t + C_1C_3 \sum_{p=0}^{2\tau} \frac{1}{\mu
(t-p+1) + u} a_{t-p} + \frac{1}{\mu (t+1) + u} C_2.
\end{split}
\end{equation}

Applying telescopic cancellation again to (\ref{eqn:thm1basic2})
from $t = 0$ to $t = T-1$, we have
\begin{equation}\label{eqn:thm1basic3-0}
\begin{split}
(\mu T + u)a_{T} \leq & u a_0 + C_1C_3 \sum_{t=0}^{T-1}
\sum_{p=0}^{2\tau} \frac{1}{\mu(t-p+1) + u} a_{t-p} + \sum_{t=0}^{T-1}\frac{1}{\mu (t+1) + u} C_2 \\
\leq & u a_0 + 2C_1C_3\tau \sum_{t=0}^{T-1} \frac{1}{\mu (t+1) +
u} a_t + \sum_{t=0}^{T-1}\frac{1}{\mu (t+1) + u} C_2.
\end{split}
\end{equation}
Substituting (\ref{ineqn:atbnd}) into \eqref{eqn:thm1basic3-0}
yields
\begin{equation}\label{eqn:thm1basic3-1}
\begin{split}
(\mu T + u)a_{T} \leq u a_0 + \frac{4}{\mu}C_1C_3\tau (a_0 + C_2
O(1)) + C_2 O(\log T),
\end{split}
\end{equation}
and consequently
\begin{equation}\label{eqn:thm1basic3}
\begin{split}
a_{T} \leq \frac{u a_0 + \frac{4}{\mu}C_1C_3\tau (a_0 + C_2 O(1))
+ C_2 O(\log T)}{\mu T + u} = O\left(\frac{\log T}{T}\right),
\end{split}
\end{equation}
which completes the proof.

\hfill\eproof

\begin{customthm}2
Suppose that the cost function of (1) satisfies the following
conditions: $f(\mathbf{x})$ is strongly convex with constant $\mu$
and $h(\mathbf{x})$ is convex; $f(\mathbf{x})$ is differentiable
and $\triangledown\! f(\mathbf{x})$ is Lipschitz continuous with
constant $L$; $\mathbb{E} \|\triangledown\!f_i(\mathbf{x}) -
\triangledown\!f(\mathbf{x})\|_2^2 \leq C_f$; $\|\partial
h(\mathbf{x})\|_2^2 \leq C_h$. Define the optimal solution of (1)
as $\mathbf{x}^*$. At time $t$, fix the step-size of the DAP-SGD
recursion (8) $\eta_t$ as $\eta = O(1/\sqrt{T})$, where $T$ is the
maximum number of iterations. Define the iterate generated by (8)
at time $t$ as $\mathbf{x}_t$. Then the running average iterate
generated by (8) at time $T$, denoted by
$$\bar{\mathbf{x}}_T = \frac{1}{T+1}\sum_{t = 0}^T \mathbf{x}_t,$$
satisfies
\begin{equation}
\begin{split}
\mathbb{E}\|\bar{\mathbf{x}}_T - \mathbf{x}^*\|_2^2 \leq
O\left(\frac{1}{\sqrt{T}}\right).
\end{split}
\end{equation}
\end{customthm}

\noindent\textbf{Proof of Theorem 2:} We start
from~(\ref{eqn:main3})  in the proof of Theorem \ref{thm:sc}.
Define the step-size rule
\begin{equation}\label{eqn:stepsizerule-new}
\begin{split}
\eta_t = \eta = \frac{1}{v\sqrt{T}},
\end{split}
\end{equation}
where $v$ is a positive constant such that $\min (\mu/(4C_4\tau),
1/L)\geq \eta$. Defining constants
$$C_4 = (2 + 56\tau)L^2,$$
and
$$C_5 = (16 \tau + 8) C_f + (43 \tau + 23) C_h,$$ followed by
manipulating (\ref{eqn:main3}), we have (similar to the inequality
\eqref{eqn:thm1basic}) the following result
\begin{equation}\label{eqn:thm2basic}
\begin{split}
a_{t+1} \leq (1 - \mu\eta) a_t + C_4 \eta^2 \sum_{p = 0}^{2\tau}
a_{t-p} + C_5 \eta^2
\end{split}
\end{equation}

Applying telescopic cancellation to (\ref{eqn:thm2basic}) from $t
= 0$ to $t = T$ yields
\begin{equation} \label{eqn:xxxx-1}
\begin{split}
a_{T+1} \leq & a_0 - \sum_{t=0}^T \mu \eta a_t + C_4 \eta^2
\sum_{t=0}^T \sum_{p = 0}^{2\tau} a_{t-p} + C_5
(T+1) \eta^2 \\
\leq & a_0 - \sum_{t=0}^T (\mu\eta - 2C_4 \eta^2\tau)a_t + C_5
(T+1) \eta^2.
\end{split}
\end{equation}
Since $\mu/(4C_4\tau) \geq \eta$ such that
\begin{equation} \label{eqn:xxxx-2}
\begin{split}
\sum_{t = 0}^T (\mu\eta - 2C_4\tau\eta^2) a_t \geq
\frac{\mu\eta}{2} \sum_{t = 0}^T a_t,
\end{split}
\end{equation}
\eqref{eqn:xxxx-1} implies
\begin{equation}\label{eqn:xxxx-3}
\begin{split}
\frac{\mu\eta}{2} \sum_{t = 0}^T a_t \leq & a_0 - a_{T+1} + C_5
(T+1) \eta^2,
\end{split}
\end{equation}
and consequently
\begin{equation}\label{eqn:thm2basic-xx}
\begin{split}
\frac{\mu\eta}{T+1} \sum_{t = 0}^T a_t \leq & \frac{2a_0 + 2C_5
(T+1) \eta^2}{T+1}.
\end{split}
\end{equation}

According to Jensen's inequality, we have
\begin{equation}\label{eqn:jensen}
\begin{split}
\frac{\mu\eta }{T+1}\sum_{t = 0}^T a_t =& \frac{\mu\eta}{T+1} \sum_{t = 0}^T \mathbb{E}\|\mathbf{x}_t - \mathbf{x}^*\|_2^2 \\
\geq & \mu\eta \mathbb{E} \left\|\frac{1}{T+1}\sum_{t= 0}^T\mathbf{x}_t - \mathbf{x}^* \right\|_2^2 \\
= & \mu\eta \mathbb{E}\|\bar{\mathbf{x}}_T - \mathbf{x}^*\|_2^2.
\end{split}
\end{equation}

Substituting (\ref{eqn:jensen}) and the step-size rule
\eqref{eqn:stepsizerule-new} into (\ref{eqn:thm2basic-xx}), we
have
\begin{equation}
\begin{split}
\mathbb{E}\|\bar{\mathbf{x}}_T - \mathbf{x}^*\|_2^2 \leq & \frac{2a_0 v\sqrt{T} + 2C_5 (T+1) \frac{1}{v\sqrt{T}}}{\mu (T+1)} = O(\frac{1}{\sqrt{T}}),
\end{split}
\end{equation}
which completes the proof. \hfill\eproof

\end{document}